\documentclass[UTF-8,reqno]{amsart}
\usepackage{enumerate}
\usepackage{tabularx,array}\usepackage{tabulary}
\usepackage{xltabular}\usepackage{xltabular}\usepackage{longtable}\usepackage{lipsum} 
\usepackage{verbatim,listings}
\usepackage{amssymb,url,color, booktabs}

\usepackage{mathrsfs}
\usepackage{amsmath}
\usepackage{verbatim}
\usepackage{fancyhdr}
\usepackage[nobysame]{amsrefs}
\BibSpec{article}{%
	+{}{\PrintAuthors} {author}
	+{,}{ \textrm} {title}
	+{.}{ \textit} {journal}
	+{,}{ \textbf} {volume}
	+{}{ \parenthesize} {date}
	+{,}{ } {pages}
	+{.}{ arXiv:} {eprint}
	+{.}{} {transition}
}
\BibSpec{book}{%
	+{}{\PrintAuthors} {author}
	+{,}{ \textit} {title}
	+{.}{ \textrm} {series} 
	+{,}{ Vol.} {volume}
	+{.}{ } {publisher}
	+{,}{ } {date}
	+{.}{} {transition}
}
\usepackage{color}
\usepackage{tikz}
\usepackage{graphicx}
\usepackage[colorlinks=true]{hyperref}
\hypersetup{
	linkcolor=blue,          
	citecolor=red,        
	filecolor=blue,      
	urlcolor=cyan
}

\numberwithin{equation}{section}

\newcommand{\be}{\begin{eqnarray}}
	\newcommand{\ee}{\end{eqnarray}}
\newcommand{\ce}{\begin{eqnarray*}}
	\newcommand{\de}{\end{eqnarray*}}
\newtheorem{theorem}{Theorem}[section]

\newtheorem{remark}{Remark}[section]

\newtheorem{lemma}[theorem]{Lemma}
\newtheorem{definition}{Definition}[section]
\newtheorem{proposition}[theorem]{Proposition}
\newtheorem{Examples}[theorem]{Example}
\newtheorem{corollary}[theorem]{Corollary}

\newenvironment{proof of theorem 1.4}{{\it Proof of Theorem 1.4}.}{{\hfill 	
		$\square$\hskip - \parfillskip}}
\newenvironment{proof of theorem 1.5}{{\it Proof of Theorem 1.5}.}{{\hfill 	
		$\square$\hskip - \parfillskip}}
\newenvironment{proof of theorem 1.6}{{\it Proof of Theorem 1.6}.}{{\hfill 	
		$\square$\hskip - \parfillskip}}

\makeatletter

\newcommand{\Rmnum}[1]{\expandafter\@slowromancap\romannumeral #1@}
\makeatother

\usepackage{caption}
\def\[{{\Big[}}
\def\]{{\Big]}}
\def\<{{\langle}}
\def\>{{\rangle}}
\def\({{\Big(}}
\def\){{\Big)}}

\def\bx{{\mathbf{x}}}

\def\min{{\mathord{{\rm min}}}}

\def\={&\!\!=\!\!&}

\def\1{{\mathbf{1}}}

\def\geq{\geqslant}
\def\leq{\leqslant}
\def\ge{\geqslant}
\def\le{\leqslant}

\def\[{{\Big[}}
\def\]{{\Big]}}
\def\<{{\langle}}
\def\>{{\rangle}}
\def\({{\Big(}}
\def\){{\Big)}}

\def\bx{{\mathbf{x}}}

\def\W{{\mathcal W}}

\def\min{{\mathord{{\rm min}}}}

\def\={&\!\!=\!\!&}
\def\bt{\begin{theorem}}
	\def\et{\end{theorem}}
\def\bl{\begin{lemma}}
	\def\el{\end{lemma}}
\def\br{\begin{remark}}
	\def\er{\end{remark}}
\def\bx{\begin{Examples}}
	\def\ex{\end{Examples}}
\def\bd{\begin{definition}}
	\def\ed{\end{definition}}
\def\bp{\begin{proposition}}
	\def\ep{\end{proposition}}
\def\bc{\begin{corollary}}
	\def\ec{\end{corollary}}

\def\geq{\geqslant}
\def\leq{\leqslant}
\def\ge{\geqslant}
\def\le{\leqslant}

\def\<{\langle} \def\>{\rangle}

\def\bpf{\begin{proof}}
	\def\epf{\end{proof}}

\allowdisplaybreaks

\tikzset{
	pattern size/.store in=\mcSize,
	pattern size = 5pt,
	pattern thickness/.store in=\mcThickness,
	pattern thickness = 0.3pt,
	pattern radius/.store in=\mcRadius,
	pattern radius = 1pt}
\makeatletter
\pgfutil@ifundefined{pgf@pattern@name@_ak5ydppfi}{
	\pgfdeclarepatternformonly[\mcThickness,\mcSize]{_ak5ydppfi}
	{\pgfqpoint{0pt}{0pt}}
	{\pgfpoint{\mcSize+\mcThickness}{\mcSize+\mcThickness}}
	{\pgfpoint{\mcSize}{\mcSize}}
	{
		\pgfsetcolor{\tikz@pattern@color}
		\pgfsetlinewidth{\mcThickness}
		\pgfpathmoveto{\pgfqpoint{0pt}{0pt}}
		\pgfpathlineto{\pgfpoint{\mcSize+\mcThickness}{\mcSize+\mcThickness}}
		\pgfusepath{stroke}
}}
\makeatother
\tikzset{every picture/.style={line width=0.75pt}} 

\begin{document}
\title{Stability of $L^p$-Brunn--Minkowski inequalities for dual 	quermassintegrals}
\author{Shanwei Ding, Jinyu Gao$^{*}$, and Guanghan Li}\thanks{* Corresponding author}


\begin{abstract}
	In this paper, we establish stability results of   $L^p$-Minkowski inequality for dual  quermassintegrals $\widetilde{W}_q$ (or dual entropy $\widetilde{E}$ when $q=0$), which in turn yields the corresponding stability of the dual $L^p$-Brunn--Minkowski inequality for  
	convex bodies $K, L \subset \mathbb{R}^{n+1}$.
	We work under the following general assumptions: the convex bodies are $C^2$-close to a ball (without assuming origin-symmetry)  and $p> q-2n-2$. For such convex bodies, we establish quantitative lower bounds for the stability deficit in terms of two distinct metrics: the $L^2$ distance and the relative asymmetry index. These estimates are strengthened versions of the corresponding inequalities. The latter, in turn, imply the uniqueness result of ($p,q$)-th dual curvature measure when $p> q-2n-2$. 
\end{abstract}

\maketitle
\setcounter{tocdepth}{2}
\tableofcontents

\section{Introduction}\label{sec 1}
The Brunn--Minkowski theory and its associated Minkowski problem are of central interest in the subject of convex geometric analysis, with deep connections to analysis, probability, and partial differential equations.
 Its development spans more than a century and has had a profound influence across mathematics. At its core lies the classical Brunn--Minkowski inequality (see \cite{book-convex-body,Gardner-book-2006,Gardner-2002,Borell-1997}) and its $L^p$ analogue, the so-called $L^p$-Brunn--Minkowski inequality introduced by \cite{Firey-book} and later developed by  \cite{Lutwak-Yang-Zhang-2012,L93}. Let $\mathcal{K}^{n+1}$ with $n\ge 1$ denote the set of convex bodies (compact, convex sets with non-empty interiors)  in $\mathbb{R}^{n+1}$, $\mathcal{K}_o^{n+1}$ denote the set of convex bodies containing the origin in their interiors.
We denote by $\mathbb{B}^{n+1}:=\{x\in\mathbb{R}^{n+1}:|x| \leq 1\}$ the standard unit ball, and by $\mathbb{S}^{n}:=\partial \mathbb{B}^{n+1}$ the unit sphere.  

For  $K\in\mathcal{K}^{n+1}$, its support function $u_K: \mathbb{S}^{n} \rightarrow \mathbb{R}$ is given by
$$
u_K(x)=\max _{y \in K}\langle x, y\rangle, \quad x \in \mathbb{S}^n,
$$
where $\langle\cdot, \cdot\rangle$ is the standard Euclidean inner product on $\mathbb{R}^{n+1}$.
It is well known that 
the support function admits a convex, positively homogeneous extension of degree one to $\mathbb{R}^{n+1}$. In addition, any convex and positively homogeneous function of degree one is the support function of a unique convex body. 
Given a function $f: \mathbb{S}^{n} \rightarrow(0, \infty)$, the set $$\left\{Q\in\mathcal{K}_o^{n+1}:u_Q\leq  f\right\}$$
has a unique maximal element,  the \textbf{Aleksandrov body} $W(f)$ associated with the function $f$ (more detial, cf. \cite{L93
}*{Page 142}). 
It is clear that if $f$ is a support function of a convex body $K$, then the $W(f)=K$. 
Let $K, L\in\mathcal{K}_o^{n+1}$.  For $p \in \mathbb{R},\, t \in(0,1)$, define the $L^p$-Minkowski convex combination (or simply $L^p$-Minkowski sum) $(1-t) \cdot K+_p t \cdot L$ of $K$ and $L$ with respect to $t$ by
$$
(1-t) \cdot K+_p t \cdot L:=\left\{\begin{array}{ll}
	W\left(\left((1-t) u_K^p+t u_L^p\right)^{\frac{1}{p}}\right), & p \neq 0, \\
	W\left(u_K^{1-t} u_L^t\right), & p=0.
\end{array} \right.
$$
Note that the case $p=1$ corresponds to the classical Minkowski sum. That is, $(1-t) \cdot K+_1 t \cdot L= (1-t) K+t L=\{(1-t) x+t y: x \in K, y \in L\}$.
By \cite{L93}*{Corollary 1.3}, the $L^p$-Brunn--Minkowski inequality states that for $p \geq 1$, $t \in (0,1)$, 
\begin{align}\label{equ:BM-neq-p=1}
	\left| (1-t) \cdot K+_p t \cdot L\right| \geq\left((1-t) \left| K\right| ^{\frac{p}{n+1}}+t \left| L\right| ^{\frac{p}{n+1}}\right)^{\frac{n+1}{p}},
\end{align}
holds with equality if and only if $K$ and $L$ are dilates (and may differ by a translation in the case $p=1$). Here $|K |$ denotes the volume of $K\subset \mathbb{R}^{n+1}$. In the case $p=1$, it reduces to the classical Brunn--Minkowski inequality. The case $p=0$ in \eqref{equ:BM-neq-p=1} is understood in the limiting sense as $p\to 0$ and is known as the log-Brunn--Minkowski inequality
\begin{align}\label{equ:BM-neq-p=0}
	\left| (1-t) \cdot K+_0t\cdot L\right| \geq \left| K\right|^{1-t}\left| L\right| ^{t}.
\end{align} 
It was proved in dimension two for origin-symmetric convex bodies in \cite{Boroczky-Lutwak-Yang-Zhang-2012}. In higher dimensions and non-origin-symmetric case, partial results have been obtained under additional assumptions in \cite{Xi-Leng-2016,Boroczky-Kalantzopoulos-2022}.

In recent years, a stronger versions of the classical Brunn--Minkowski inequality, namely, the case $p<1$ for $L^p$-Brunn--Minkowski inequality and its equivalent form, $L^p$-Minkowski inequality (see \cite{Boroczky-Lutwak-Yang-Zhang-2012,Milman-JEMS-2025}), 
\begin{align}
	\label{equ:M-neq}
	\frac{1}{p}\int_{\mathbb{S}^n}u_L^pu_K^{1-p}\sigma_{n}(\lambda_K)~dx\geq\frac{n+1}{p} \left| K\right| ^{\frac{n+1-p}{n+1}}\left| L\right| ^{\frac{p}{n+1}},\quad &\mathrm{ if }\ p\ne0,\\
	\label{equ:M-neq=p=1}
	\int_{\mathbb{S}^n}\left(\log(\frac{u_L}{u_K})\right)u_K\sigma_n(\lambda_K) ~dx
	\geq  \left| K\right| \log\frac{\left| L\right| }{\left| K\right| },\quad \ \! &\mathrm{ if }\ p=0,
\end{align} 
has attracted considerable attention; see, for example, \cite{Boroczky-Lutwak-Yang-Zhang-2012,Chen-Huang-Li-Liu-2020,Putterman-2021,Milman-JEMS-2025,Kolesnikov-Milman-MAMS-2022,Milman-JDG-2024,Boroczky-Lutwak-Yang-Zhang-2013,Yang-Zhang-2019}. 
Here, $\lambda_K$ denotes the principal radii of curvature of the hypersurface $\partial K$, and $\sigma_n(\lambda_K)$ is the product of all principal radii of curvature. If $K$ is convex, then $\lambda_K$ are the eigenvalues of the matrix $\tau[u_K]=\nabla^2u_K+u_K I$, and $\sigma_n(\lambda_K)=\sigma_n(\tau[u_K])=\det(\tau[u_K])$. Moreover, $dx=d\sigma_{\mathbb{S}^n}(x)$ denotes the standard (induced) volume element of the sphere $\mathbb{S}^n$, and $\nabla$ denotes the Levi-Civita connection of the induced metric on $\mathbb{S}^n$.
For $0<p<1$, as shown by counterexamples in \cite{Boroczky-Lutwak-Yang-Zhang-2012}, the $L^p$-Brunn--Minkowski inequality and the $L^p$-Minkowski inequality do not hold in general. However, it was conjectured and later proved partially in \cite{Boroczky-Lutwak-Yang-Zhang-2012,Kolesnikov-Milman-MAMS-2022,Chen-Huang-Li-Liu-2020,Putterman-2021} that it holds under the additional assumption that $K$ and $L$ are origin-symmetric.


For $p<0$, the problem is difficult even for origin-symmetric bodies. 
Known counterexamples constructed in \cite{Boroczky-Saroglou-Cvpde-2024,Kolesnikov-Milman-MAMS-2022,Milman-JDG-2024} demonstrate that the $L^p$-Brunn--Minkowski and $L^p$-Minkowski inequalities fail for some origin-symmetric bodies when $p<0$. Moreover, It was proved in \cite{Kolesnikov-Milman-MAMS-2022,Boroczky-Saroglou-Cvpde-2024} that uniqueness fails in the even $L^p$-Minkowski problem for certain $p<0$, where the $L^p$-Minkowski problem asks for a convex body whose $L^p$-surface area measure equals a prescribed measure, and “even” means that 
the solutions are restricted to origin-symmetric convex bodies. In particular, Böröczky and Saroglou \cite{Boroczky-Saroglou-Cvpde-2024} proved that, 
 the set of non-unique solutions to the $L^p$
Minkowski problem is dense in the class of convex bodies containing the origin if $-(n+1)<p<0$.





The $L^p$-Brunn--Minkowski inequality has inspired many generalizations, often referred to as Brunn--Minkowski-type inequalities, in which volume is replaced by other functions. For instance, partial results have been obtained for the intrinsic volume  by \cite{Bianchini-Colesanti-Pagnini-Roncoroni-Math ann-2023, Patsalos-Saroglou-arxiv,Ding-Gao-Li-BM1},
and for the dual quermassintegral $\widetilde{W}_q$ (defined in \eqref{equ:dual quermassintegral}) by \cite{Sadovsky-ZhangGaoYong-2025,Xi-Zhang-2022}. 

In 1975, Lutwak \cite{Lutwak-1975} introduced the dual mixed volume, which in a special case reduces to the $j$-th dual quermassintegral, $j \in {1,\dots,n}$, and author also contributed the dual Burnn-Minkwski inequality by replacing the classical sum by radial sum. Subsequently, for a general $q \in \mathbb{R} \setminus {0}$, the \textbf{$q$-th dual quermassintegral} \cite{Sadovsky-ZhangGaoYong-2025}*{eq. (3)} was defined as follows:
\begin{align}\label{equ:dual quermassintegral}
	\widetilde{W}_q(K)=\int_{\mathbb{S}^n}\rho_K^q(\xi) d\xi,
\end{align}
where $\rho_K(\xi)=\max\{\lambda:\lambda \xi\in K\}$ with  $\xi\in \mathbb{S}^{n}$ is  the radial function of $K$ and $d\xi:=d\sigma_{\mathbb{S}^n}(\xi)$. Specially, when $q=n+1$, $\widetilde{W}_{n+1}(K)=(n+1)|K|$.
If $K$ is of class $C^2$, let $\mathscr{A}^{*} : \mathbb{S}^n \to \mathbb{S}^n$ be the reverse radial Gauss map of $K$ (for details, see \cite[eq. (1.5)]{Chen-Li-2021}). Since the determinant of the Jacobian of $\mathscr{A}^{*}$ is given by \cite[eq. (2.28)]{Ding-Li-Trans}  
\begin{align}\label{equ:dx=Adxi}
	|\operatorname{Jac}\mathscr{A}^{*}|(x)=\left| \frac{d\xi}{dx}\right| =\frac{u_K(x)\sigma_n(\lambda_K(x))}{\rho_K^{n+1}\left(\mathscr{A}^{*}(x)\right)},
\end{align}
then $\widetilde{W}_q(K)$ can be expressed in the following form:
\begin{align}\label{equ:dual quermassintegral-C2}
	\widetilde{W}_q(K)=\int_{\mathbb{S}^n}u_K(x)\rho_K^{q-n-1}(\mathscr{A}^*(x))\sigma_n(\lambda_K(x))~dx,\quad q\neq 0.
\end{align}
The case $q=0$ is singular, and   the corresponding \textbf{dual entropy} \cite{Sadovsky-ZhangGaoYong-2025}*{Page 3} is defined as 
\begin{align}\label{equ: dual entropy}
	\widetilde{E}(K)=\int_{\mathbb{S}^n}\log \rho_K(\xi) d\xi,
\end{align}
If $K$ is of class $C^2$, then $\widetilde{E}(K)$ can be expressed in the following form:
\begin{align}\label{equ: dual entropy-C2}
	\widetilde{E}(K)=\int_{\mathbb{S}^n}(\log \rho_K) u_K\rho_K^{-n-1}\sigma_n(\lambda_K)~dx.
\end{align}

In \cite{Huang-Lutwak-Yang-Zhang-2016}, Huang, Lutwak, Yang, and Zhang proposed the dual Minkowski problem, and in \cite{Lutwak-Yang-Zhang-2018}, Lutwak, Yang, and Zhang introduced the $(p,q)$-th dual curvature measure $\tilde{C}_{p,q}(K,\cdot)$.   
There is a close relationship among Minkowski type inequalities, Brunn--Minkowski type inequalities, and the Minkowski type problem. For $q=n+1$, $p>-(n+1)$, and $K,L$ are origin-symmetric convex bodies, the fact that the $L^p$-Brunn--Minkowski inequality, the uniqueness of solutions to the $L^p$-Minkowski problem, and the $L^p$-Minkowski inequality  are equivalent has been proved in  \cite{Boroczky-Lutwak-Yang-Zhang-2012,Milman-JEMS-2025}.

In \cite{Xi-Zhang-2022}, Xi and Zhang established dual Brunn--Minkowski inequalities for all real numbers $p \ge q$, based on the dual $L_p$-Minkowski inequalities they had obtained. Using these dual Minkowski inequalities, they further derived the uniqueness of the $(p,q)$-dual curvature measure. Earlier, the uniqueness of solutions for the case $p > q$ had already been proved by employing the strong maximum principle; see, e.g., Huang and Zhao \cite[Theorem 1.5]{Huang-Zhao-2018}. For the case of polytopes with 
$p \geq q$, the uniqueness was established by Lutwak, Yang, and Zhang \cite[Theorem 8.3]{Lutwak-Yang-Zhang-2018}. And \cite{Wang-Zhou-2026} further extend the important result in \cite{Lutwak-Yang-Zhang-2018}  to the setting of general convex bodies.
\begin{theorem}[\cite{Xi-Zhang-2022}*{Theorem A}]
	Let $K, L\in\mathcal{K}^{n+1}_o$. Suppose $p \geq q$ are real numbers. Then  	
	\begin{align}
		\label{equ:dual BM pq}
		\widetilde{W}_{q}\left((1-t) \cdot K+_p t \cdot L\right)^{\frac{1}{q}} \geq\left((1-t) \widetilde{W}_{q}(K)^{\frac{p}{q}}+t \widetilde{W}_{q}(L)^{\frac{p}{q}}\right)^{\frac{1}{p}},
	\end{align}
	for $q \neq 0$ and $t \in(0,1)$.
	Equality holds if and only if $K$ and $L$ are dilates. When $p=0$, the right side will be understood as the limit, and  	  	
	\begin{align}
		\label{equ:dual BM 0q}
		\widetilde{W}_{q}\left((1-t) \cdot K+_0 t \cdot L\right) \leq \widetilde{W}_{q}(K)^{1-t} \widetilde{W}_{q}(L)^t, \quad \text { for } q \neq 0 \text { and } t \in(0,1),
	\end{align}
	with equality if and only if $K$ and $L$ are dilates.
	When $q=0$, we have  	
	\begin{align}\label{equ:dual BM p0}
		\widetilde{E}\left((1-t) \cdot K+_p t \cdot L\right) \geq(1-t) \widetilde{E}(K)+t \widetilde{E}(L).
	\end{align}
	If $p>0$, equality holds if and only if $K=L$; if $p=0$, equality holds if and only if $K$ and $L$ are dilates.
\end{theorem}
When we take $q = n+1$ in \eqref{equ:dual BM pq} and \eqref{equ:dual BM 0q}, they reduce to the classical $L^p$-Brunn--Minkowski inequality and the log-Brunn--Minkowski inequality, i.e., \eqref{equ:BM-neq-p=1} and \eqref{equ:BM-neq-p=0}, respectively.

The case  where $q\neq n+1$ and $p<q$ is considerably more difficult. A significant advance was made by Sadovsky and Zhang \cite{Sadovsky-ZhangGaoYong-2025}, who proved the dual Brunn--Minkowski inequalities for origin-symmetric convex bodies $K,L$ in the range $p=1$, $0<q\le n+1$, thereby confirming a conjecture of Huang, Lutwak, Yang, and Zhang (as noted in \cite{Huang-Yang-Zhang-2025}). 
Following Sadovsky and Zhang \cite{Sadovsky-ZhangGaoYong-2025} and applying H\"older's inequality, one obtains the dual Brunn--Minkowski inequalities for $p > 1$, $0 < q \le n+1$ in the origin-symmetric case (see \cite[Theorem 1.4 (i)]{Chen-Tang-Wang-2026-arxiv}). Furthermore, if $K$ satisfies $\|u_K - 1\|_{C^{2,\alpha}} \le \varepsilon$, then Theorem 1.4 (ii) of \cite{Chen-Tang-Wang-2026-arxiv} gives the corresponding result for $0 < p < 1$, $p < q \le n+1$ in the origin-symmetric case.
Hu and Ivaki \cite{Hu-Ivaki-2026-arxiv} considered the unconditional setting (a stronger symmetry condition than origin-symmetry), where a convex body is symmetric with respect to every coordinate hyperplane. They established in \cite[Theorem 1.3]{Hu-Ivaki-2026-arxiv} the dual Brunn--Minkowski inequality with $p=0$, $q\in(0,n+1)$, and extended it in \cite[Theorem 1.4]{Hu-Ivaki-2026-arxiv} to $0<p<1$, $q=n+1+\alpha$ with $0<\alpha\leq \frac{2p(1-p)}{2-p}$. And for case $p=1$, in \cite{Li-Wan-Zhang-arxiv-2026}, Li, Wan, and Zhang further investigated the  Brunn--Minkowski inequality for $q$-th dual quermassintegrals. In particular, they established the inequality for unconditional convex bodies in the full range $0<q\leq n+2$. They also proved the endpoint case $q=n+3$ for origin-symmetric convex bodies, while showing that the inequality fails for general convex bodies when $q>n+1$, and even in the origin-symmetric class when $q>n+3$. As applications, they derived several uniqueness results for the corresponding dual curvature measures.
 In \cite{Tao-2026}, Tao determined the sharp $L_p$ Brunn--Minkowski range for planar dual quermassintegrals on origin-symmetric convex bodies: $0 < q \leq 2$ for $0 \leq p < 1$, and $0 < q \leq 4$ for $p = 1$. 
 In \cite{Xiong-Yang-2026}, Xiong and Yang established the dual Brunn--Minkowski inequality for $p=0$ and origin-symmetric convex bodies in the case $n=1$, $0<q<n+1$, and showed that it fails for $q>n+1\ge 2$.


Besides, Chen, Chen, and Li \cite[Theorem 1.5]{Chen-Chen-LiQiRui-2021} provides a counterexample showing that  when $p<q-2n-2$, even  for $K = \mathbb{B}^{n+1}$ and an origin-symmetric convex body $L$ that is $C^0$-close to $\mathbb{B}^{n+1}$, the opposite inequality holds in \eqref{equ:dual BM pq}, which equivalently serves as a counterexample to the uniqueness of the corresponding $(p,q)$-th   dual Minkowski problem \cite{Chen-Chen-LiQiRui-2021}*{Theorem 1.1}.  
Consequently, it is appropriate to restrict to the case $p \geq q - 2n - 2$ when considering the dual $L^p$-Brunn--Minkowski inequality and the uniqueness of the corresponding dual $L^p$-Minkowski Problem. 

In the isotropic case (i.e., $L=\mathbb{S}^n$), Chen, Huang, and Zhao \cite[Theorem 1.3]{Chen-Huang-Zhao-2019} proved the uniqueness of the  $(p,q)$-th   dual Minkowski problem for $p\geq -n-1$, $q \leq \min\{n+1,\, n+1+p\}$, and $p \neq q$, where $K$ is an origin-symmetric convex body. Also in the isotropic setting, Chen and Li \cite[Theorem 1.4 and 7.1]{Chen-Li-2021} proved uniqueness for $p>1$ with $q\leq n+1$, and for $-n-1\leq p< q<1$.  Ivaki and Milman \cite[Corollary 1.5]{Ivaki-Milman-2023} established uniqueness of $(p,q)$-th   dual  Minkowski problem when $-n-1\leq p$ and $q\leq n+1$, with at least one of the inequalities being strict.


Li, Sheng, and Wang \cite[Theorem 1.6]{Li-Sheng-Wang-2022} used a flow method to discuss the uniqueness for the anisotropic and isotropic
$(p,q)$-th even  dual   Minkowski problem when $p=1$ and various $q$, while Ding and Li \cite{Ding-Li-Trans}*{Corollary 1.7 and Theorem 1.8} used the convergence of a flow to derive the uniqueness results for the  anisotropic and isotropic $(p,q)$-th    dual Minkowski problem with various $p,q$. 
In \cite{Li-Wan-JFA-2024}, Li and Wan employed various methods to study the uniqueness of $(p,q)$-th   dual Minkowski problem for $-n-1< p\leq -1$ and $n+1 \leq q \leq n+\frac{1}{2}+ \sqrt{\frac{1}{4}-\frac{(1+p)(n+1+p)}{n(n+2)}}$, subsequently, in \cite{Li-Wan-JRAM-2026}, they completed the classification of solutions to the planar isotropic $(p,q)$-th   dual Minkowski problem.

For general $L$ (the anisotropic case), \cite{Li-Liu-Lu-2022} established the nonuniqueness of solutions to the $(p,q)$-th   dual Minkowski problem for $p < 0 < q$. Uniqueness results in the anisotropic setting have also been obtained by Cabezas-Moreno and Hu \cite{Cabezas-Hu-2025} for $(p,q)$-th  even  dual Minkowski problem with $1 < p < q \leq n+1$, Hu and Ivaki \cite{Hu-Ivaki-2024},  B\"{o}r\"{o}czky, Chen,  Liu, and Saroglou \cite{Boroczky-Chen-Liu-Saroglou-arxiv} for ``near isotropic" $(p,q)$-th  even  dual Minkowski problem for  $-1<p<q<\min\{n+1,n+1+p\}$
and $q>0$, also  \cite{Boroczky-Chen-Liu-Saroglou-arxiv-2}*{Corollary 1.3} for $\mathbb{S}^2$ and $p\in(0,1)$, $|q-3|<\varepsilon$. 
Given the above examples and the literature, it is challenging to drop the even measure condition (namely, the condition that the convex bodies are origin‑symmetric) when $p<q$.

In this paper, we investigate the dual Brunn--Minkowski conjecture proposed by Lutwak et al. from three perspectives. First, we extend the admissible range of the parameters 
$p$ and $q$ beyond the traditional restriction 
$p\geq q$ to case $p > q - 2n - 2$. Second, we relax the origin‑symmetry assumption on the convex bodies 
$K$ and 
$L$; in particular, the case 
$p<q$ remains highly challenging without this condition, and counterexamples are known to exist. Third, we establish stability (i.e., sharp quantitative) estimates for the inequality.


Inspired by \cite{Boroczky-Chen-Liu-Saroglou-arxiv,Boroczky-Saroglou-Cvpde-2024,Patsalos-Saroglou-arxiv,Bianchini-Colesanti-Pagnini-Roncoroni-Math ann-2023,Ding-Gao-Li-BM1,Chen-Chen-LiQiRui-2021}, we consider the convex bodies that are $C^2$-close to the unit ball. 
We denote by $\mathcal{N}_{\varepsilon}$ the collection of all $C^2$ convex bodies $K$ satisfying
$$
\| u_K - 1 \|_{C^2(\mathbb{S}^n)} < \varepsilon.
$$
Here  for $m\in\mathbb{Z}^+$, we denote   $$\| f \|_{C^m}:=\| f \|_{C^m(\mathbb{S}^n)}=\max_{0\leq k\leq m}\sup_{x\in\mathbb{S}^n}\left|\nabla^k f(x)\right|.$$

In this context, the stability of an inequality concerns how close the sets $K$ and $L$ are to being homothetic when the equality is nearly attained. To quantify this closeness,
we introduce the following two distance measures to quantify the deviation between $K$ and $L$:
(i) the classical 
$L^2$ distance \eqref{equ:L2-distance} (cf. \cite{book-convex-body}*{eq. (7.120)}) \begin{align}\label{equ:L2-distance}
	\delta_2(K,L)=\left\Vert u_K-u_L\right\Vert_{L^2}=\left(\int_{\mathbb{S}^n}|u_K-u_L|^2~dx\right)^{\frac{1}{2}},
\end{align} 
(ii) the more natural \textbf{relative asymmetry index} used by Figalli, Maggi, and Pratelli \cite{Figalli-Maggi-Pratelli-2010}*{eq.  (1.29)}, which defined by 
\begin{align}
	\label{equ:def-asymmetry}
	\alpha(K,L)=\min _{y \in \mathbb{R}^{n+1}}\left\{\frac{|K \Delta(r L+y)|}{|K|}:|r L|=|K|\right\},
\end{align}
where  $K \Delta L=(K \backslash L) \cup(L \backslash K)$ denotes the symmetric difference between two sets $K,L\subset\mathbb{R}^{n+1}$. 
Distance $\alpha(K,L)$ satisfies symmetry and the triangle inequality \cite{Figalli-Maggi-Pratelli-2010}*{eq. (4.2)}.
For the case $u_L\equiv1$, i.e., $L=\mathbb{B}^{n+1}$, the relative asymmetry index is known as the Fraenkel asymmetry (ref. \cite{Fusco-Maggi-2008}).

The stability of Brunn--Minkowski type inequalities is a well-known problem that has attracted considerable attention in recent years. 
The first stability forms of the Brunn--Minkowski inequality for convex
bodies were due to Minkowski himself  and more details can be seen in \cite{Groemer-1993}. Subsequently, nearly optimal stability results were obtained by Diskant \cite{Diskant-1973} and Groemer \cite{Groemer-1988}. Later developments introduced natural distance notions and established the optimal exponent, with further extensions to general bounded measurable sets in \cite{Figalli-Maggi-Pratelli-2010,Figalli-Maggi-Pratelli-2009,Figalli-Hintum-Tiba-arxiv-2023,Figalli-Hintum-Tiba-arxiv-2024}.
Böröczky, Machado, and Ramos \cite[Lemma 4.1]{Boroczky-Machado-Ramos-2026-arxiv} proved a stability result for the $L^p$-Brunn--Minkowski inequality for  $p>1$, and for convex bodies containing the origin, measured in terms of the relative asymmetry index,  using the quantitative stability version of the Brunn--Minkowski inequality (i.e. $p=1$) for convex sets in \cite{Boroczky-Figalli-Ramos-2025-book}. 
Böröczky and De \cite{Boroczky-De-2024} proved a stability result for the $L^p$-Brunn--Minkowski inequality with $p=0$ (i.e., the logarithmic case),  assuming that the bodies are symmetric with respect to $n+1$ independent linear hyperplanes.  
To the best of our knowledge, stability results for the dual $L^p$-Brunn--Minkowski inequality in the case $q \neq n+1$ are not yet available in the literature.

We consider a more general setting where the convex bodies 
$K,L$ are not required to be origin-symmetric. To make the analysis tractable in this broader setting, we introduce the \textbf{$p$-weighted Steiner point} of a convex body $K$ as
\begin{equation}
	\label{equ:def-p-Steiner}
	S_p(K)=
	\left\lbrace 
	\begin{array}{ll}
		\int_{\mathbb{S}^n}u_K^p(x)\cdot x~dx, &\text{\ for\ } p\neq 0,
		\\
		
		\int_{\mathbb{S}^n}\log u_K(x)\cdot x~dx, &\text{\ for\ } p= 0.
	\end{array}
	\right.
\end{equation}

In this paper, we assume that $p$-weighted Steiner points of $K,L$ are   at the origin ${o}$, i.e., $S_p(K)=S_p(L)={o}$.

\begin{remark}
	(1) When $p=1$, \eqref{equ:def-p-Steiner} reduces to the classical definition of the Steiner point \cite{book-convex-body}, that is $s(K)=\frac{n+1}{|\mathbb{S}^n|}S_1(K)$. When $p=-n-2$, \eqref{equ:def-p-Steiner} reduces to the  definition of the Santalo point.
	
	(2) The condition $S_p(K)={o}$ includes the origin-symmetric case as a special case.
	
	(3) This definition of $S_p(K)$ follows the weighted centering conditions for
	the entropy point in
	\cite{Andrews-Guan-Ni-2016}*{Lemma 2.5, eq. (2.7)} and
	\cite{Guan-Ni-2017}*{Lemma 2.2, eq. (2.2)}. We use it here to handle the
	general case without assuming origin symmetry.

	(4) Even if the Steiner points of $K,L\in\mathcal{N}_{\varepsilon}$ lie at the origin, this does not in general imply that the Steiner point of their $L^p$-Minkowski convex combination
	$Q_t:=(1-t)\cdot K+_pt\cdot L\in\mathcal{N}_{O(\varepsilon)}$ with $t\in[0,1]$ is also at the origin, except in the special case $p=1$.
	In contrast, a key feature of $S_p(K)=S_p(L)={o}$ is that it is preserved under $L^p$-Minkowski convex combinations. That is, if $K$ and $L$  satisfy $S_p(K)=S_p(L)={o}$, then $S_p(Q_t)={o}$. Lemma \ref{lemma:C2-1-p-sum-C2-1} makes the above verification straightforward.
	

	(5) The condition $S_p(K)={o}$ is not preserved under translations. This is consistent with the fact that translations do not appear in the equality case of the $L^1$-Minkowski inequality stated in Theorem~\ref{thm:p-M-inequ} with $p=1,q=n+1$. 
\end{remark}

Under the condition $S_p(K)=S_p(L)={o}$, we establish stability for the dual $L^p$-Brunn--Minkowski inequality over a broader range of $p,q\in\mathbb{R}$ (i.e., $p>q-2n-2$) and for non-origin-symmetric convex bodies, with respect to several distances. The following stability result of $L^p$-Brunn--Minkowski inequalities for dual quermassintegrals can be derived from  corresponding stability results for the dual $L^p$-Minkowski inequality in Corollary \ref{thm:dual p-M-L2}.

\begin{theorem}
	\label{thm:dual p-BM-L2-Vk!=Vk}
	Suppose $K, L \in \mathcal{K}^{n+1}$ are $C^3$ convex bodies, $p$ and $q$ are real numbers with  $p > q - 2n - 2$, and $t \in [0,1]$.  For all $\eta>0$, there exists $\varepsilon>0$  such that if $K,L\in\mathcal{N}_{\varepsilon}$ and $K,L$ satisfy  $S_p(K)=S_p(L)={o}$, then for $p\neq 0$ and $q\neq0$, we have
	\begin{align}
		\label{equ:stab L2-BM-Wq!=Wq}
		&\frac{1}{p}\left(\widetilde{W}_q\((1-t)\cdot K +_p t\cdot L\)^\frac{p}{q}
		-\left((1-t)\widetilde{W}_q(K)^{\frac{p}{q}}+t\widetilde{W}_q(L)^{\frac{p}{q}}\right)\right)\nonumber
		\\
		\geq &\frac{t(1-t)\widetilde{W}_q(K)^{\frac{p}{q}}\widetilde{W}_q(L)^{\frac{p}{q}}}{(1-t)\widetilde{W}_q(K)^{\frac{p}{q}}+t\widetilde{W}_q(L)^{\frac{p}{q}}}|\mathbb{S}^n|^{-1}\left(\frac{2n+2+p-q}{2}-\eta\right)\delta_2^2(N_K K,N_L L),
	\end{align}
	and for $p=0, q\neq 0$, we have 
	\begin{align}
		\label{equ:stab L2-BM-Wq!=Wq-p=0}
		&\frac{1}{q}\left(\log\widetilde{W}_q\((1-t)\cdot K +_p t\cdot L\)-\left((1-t)\log\widetilde{W}_q(K)+t\log\widetilde{W}_q(L)\right)\right)\nonumber
		\\
		\geq &t(1-t)|\mathbb{S}^n|^{-1}\left(\frac{2n+2-q}{2}-\eta\right)\delta_2^2(N_K K,N_L L),
	\end{align}
	and for $p\neq 0,q=0$, we have
	\begin{align}\label{equ:stab L2-BM-q=0}
		&\frac{1}{p}\left(
		e^{p\frac{\widetilde{E}\left((1-t)\cdot K+_pt\cdot L\right)}{|\mathbb{S}^n|}}
		-
		\left((1-t)e^{p\frac{\widetilde{E}(K)}{|\mathbb{S}^n|}}+te^{p\frac{\widetilde{E}(L)}{|\mathbb{S}^n|}}\right)
		\right)
		\nonumber
		\\
		\geq& t(1-t)C_1
		\left(\frac{2n+2+p}{2}-\eta\right)\delta_2^2(\widetilde{N_K}K,\widetilde{N_L}L),
	\end{align}
	and for $p=q=0$, we have
	\begin{align}\label{equ:stab L2-BM-p=q=0}
		&\widetilde{E}\((1-t)\cdot K +_p t\cdot L\)-\left((1-t)\widetilde{E}(K)+t\widetilde{E}(L)\right)\nonumber
		\\
		\geq&t(1-t)\left((n+1)-\eta\right)\delta_2^2(\widetilde{N_K} K,\widetilde{N_L} L).
	\end{align}
	Here $N_K=\left(\frac{|\mathbb{S}^n|}{\widetilde{W}_q(K)}\right)^{\frac{1}{q}}, N_L=\left(\frac{|\mathbb{S}^n|}{\widetilde{W}_q(L)}\right)^{\frac{1}{q}}$, $\widetilde{N_K}=e^{-\frac{\widetilde{E}(K)}{|\mathbb{S}^n|}}, \widetilde{N_L}=e^{-\frac{\widetilde{E}(L)}{|\mathbb{S}^n|}}$,  $|\mathbb{S}^n|$ denotes the surface area of $\mathbb{S}^n$, and 
	$$C_1=\frac{e^{\frac{p\widetilde{E}(K)}{|\mathbb{S}^n|}}e^{\frac{p\widetilde{E}(L)}{|\mathbb{S}^n|}} }{\left((1-t)e^{p\frac{\widetilde{E}(K)}{|\mathbb{S}^n|}}+te^{p\frac{\widetilde{E}(L)}{|\mathbb{S}^n|}}\right)\cdot|\mathbb{S}^n|}
	.$$                 
\end{theorem} 

\begin{remark}\label{rk:L2-replaced by the anisotropic weighted-BM}
		In the case $p\neq 0$, $q=0$, if $p>0$, then by Young’s inequality, inequality \eqref{equ:stab L2-BM-q=0} becomes
	\begin{align}
		&\frac{1}{p}\left(
		e^{p\frac{\widetilde{E}\left((1-t)\cdot K+_pt\cdot L\right)}{|\mathbb{S}^n|}}
		-
	e^{p\frac{(1-t)\widetilde{E}(K)+t\widetilde{E}(L)}{|\mathbb{S}^n|}}\right)	
		\nonumber
		\\
		\geq& t(1-t)C_1
		\left(\frac{2n+2+p}{2}-\eta\right)\delta_2^2(\tilde{K},\tilde{L}),
	\end{align}
which is sharper than inequality \eqref{equ:dual BM p0},  likewise, \eqref{equ:stab a-BM-p!=q=0} in Theorem \ref{thm:dual p-BM-alpha-Vk!=Vk} is sharper than  inequality \eqref{equ:dual BM p0} for $p>0$.
	
	
\end{remark}
Theorem  \ref{thm:dual p-BM-alpha-Vk!=Vk}  can be derived  from corresponding stability results for the dual $L^p$-Minkowski inequality in Corollary \ref{thm:dual p-M-al}.


\begin{theorem}
	\label{thm:dual p-BM-alpha-Vk!=Vk}
	Suppose $K, L \in \mathcal{K}^{n+1}$ are $C^3$ convex bodies, $p$ and $q$ are real numbers with  $p > q - 2n - 2$, and $t \in [0,1]$.  For all $\eta>0$, there exists $\varepsilon>0$  such that if $K,L\in\mathcal{N}_{\varepsilon}$ and $K,L$ satisfy  $S_p(K)=S_p(L)={o}$, then for $p\neq 0$ and $q\neq0$, we have
	\begin{align}
		\label{equ:stab a-BM-Wq!=Wq}
		&\frac{1}{p}\left(\widetilde{W}_q\((1-t)\cdot K +_p t\cdot L\)^\frac{p}{q}
		-\left((1-t)\widetilde{W}_q(K)^{\frac{p}{q}}+t\widetilde{W}_q(L)^{\frac{p}{q}}\right)\right)\nonumber
		\\
		\geq &\frac{t(1-t)\widetilde{W}_q(K)^{\frac{p}{q}}\widetilde{W}_q(L)^{\frac{p}{q}}}{(1-t)\widetilde{W}_q(K)^{\frac{p}{q}}+t\widetilde{W}_q(L)^{\frac{p}{q}}}\left(c(n,p,q)-\eta\right)\alpha(N_K K,N_L L)^2,
	\end{align}
	and for $p=0, q\neq 0$, we have 
	\begin{align}
		\label{equ:stab a-BM-Wq!=Wq-p=0}
		&\frac{1}{q}\left(\log\widetilde{W}_q\((1-t)\cdot K +_p t\cdot L\)-\left((1-t)\log\widetilde{W}_q(K)+t\log\widetilde{W}_q(L)\right)\right)\nonumber
		\\
		\geq &t(1-t)\left(c(n,0,q)-\eta\right)\alpha(N_K K,N_L L)^2,
	\end{align}
	and for $p\neq 0,q=0$, we have
	\begin{align}\label{equ:stab a-BM-p!=q=0}
		&\frac{1}{p}\left(
		e^{p\frac{\widetilde{E}\left((1-t)\cdot K+_pt\cdot L\right)}{|\mathbb{S}^n|}}
		-
		\left((1-t)e^{p\frac{\widetilde{E}(K)}{|\mathbb{S}^n|}}+te^{p\frac{\widetilde{E}(L)}{|\mathbb{S}^n|}}\right)
		\right)
		\nonumber
		\\
		\geq&t(1-t)C_1|\mathbb{S}^n|\left(c(n,p,0)-\eta\right)\alpha(\widetilde{N_K} K,\widetilde{N_L} L)^2,
	\end{align}
	and for $p=q=0$, we have
	\begin{align}\label{equ:stab a-BM-p=q=0}
		&\widetilde{E}\((1-t)\cdot K +_p t\cdot L\)-\left((1-t)\widetilde{E}(K)+t\widetilde{E}(L)\right)\nonumber
		\\
		\geq&t(1-t)|\mathbb{S}^n|\left(c(n,0,0)-\eta\right)\alpha(\widetilde{N_K} K,\widetilde{N_L} L)^2.
	\end{align}
	Here 
	\begin{align*}
		c(n,p,q)=
			\frac{2n+2+p-q}{2(n+1)^2}.
	\end{align*}
\end{theorem}

The dual $L^p$-Minkowski inequality  implies an analogous dual $L^p$-Brunn--Minkowski inequality. It therefore follows that stability for the latter can be deduced from stability for the former. Accordingly, to obtain stability for the dual $L^p$-Brunn--Minkowski inequality (i.e. Theorem \ref{thm:dual p-BM-L2-Vk!=Vk} and \ref{thm:dual p-BM-alpha-Vk!=Vk}), we shall establish the  corresponding stability results for the dual $L^p$-Minkowski inequality in Section \ref{sec 4}.

Besides, the dual $L^p$-Minkowski inequality  can derive the uniqueness results for the $(p,q)$ dual curvature measures (up to a diation when $p=q$) as follows.
\begin{corollary}\label{thm:unique}
	Suppose $p$ and $q$ are real numbers with  $p > q - 2n - 2$. There exists $\varepsilon>0$, that depends only on $n,\, p$ and $q$, such that 
	if $K,L\in\mathcal{N}_{\varepsilon}$ are $C^3$ convex bodies satisfy $S_p(K)=S_p(L)={o}$,
	then the equation	
	\begin{align}\label{equ:Christolell-equ}
		\rho_K^{q-n-1}(\mathscr{A}_K^*(x))u_{K}^{1-p}(x) \sigma_n\left(\lambda_K(x\right)=\rho_L^{q-n-1}(\mathscr{A}_L^*(x))u_{L}^{1-p}(x) \sigma_n\left(\lambda_L(x)\right)
	\end{align}
	implies $K=L$ when $p\neq q$ and   $K=NL$ when $p= q$.
\end{corollary}
The case $q=n+1$ in Corollary \ref{thm:unique} corresponds to the $L^p$-Minkowski problem, while the case $q=0$ corresponds to the  $L^p$-integral curvature. 

The rest of this paper is organized as follows. 
In Section \ref{sec 3}, we study  some characteristics and formulations for  convex bodies that are $C^2$-close to a ball. 
In Section \ref{sec 4}, we   establish the stability of dual $L^p$-Minkowski inequality (i.e., Corollaries \ref{thm:dual p-M-L2} and \ref{thm:dual p-M-al}). 
Section \ref{sec 5} concerns the stability of the dual $L^p$-Brunn--Minkowski inequality (see Theorems \ref{thm:dual p-BM-L2-Vk!=Vk}  and \ref{thm:dual p-BM-alpha-Vk!=Vk}, which follows from Corollaries \ref{thm:dual p-M-L2} and \ref{thm:dual p-M-al}. Finally, in Section \ref{sec 6}, we obtain a the uniqueness results for the $(p,q)$ dual curvature measures (see Corollary \ref{thm:unique}).

\section{Some propositions for convex bodies that are $C^2$-close to a ball}\label{sec 3}

In this section, we consider  convex bodies that are $C^2$-close to a ball, and compute the corresponding geometric quantities.

For convenience, we set $w_i := \nabla_{e_i} w$ and $w_{ij} := \nabla_{e_i}\nabla_{e_j} w$ for any $w\in C^2(\mathbb{S}^n)$, with $\{e_i\}_{i=1}^n$ an orthonormal frame on $\mathbb{S}^n$. 

Now we consider convex bodies $K,L\in\mathcal{N}_{\varepsilon}$. 
We denote the support functions of $M=\partial K$ and $\W=\partial L$ are $u_K=1+w$ and $u_L=1+\bar{w}$,   respectively. We adopt the notations
\begin{align*}
	v:=u_L-u_K,
\end{align*}
and
$$\mathscr{E}:=
O(\varepsilon)\cdot \(\Vert v\Vert_{L^2}^2+\Vert \nabla v\Vert_{L^2}^2
\)
.$$


Then we can derive following lemmas
\begin{lemma}\label{lemma:rho Wq=}
	Let $q\neq 0,\alpha\in\mathbb{R}$ and  $K\in\mathcal{N}_{\varepsilon}$, then we have 
	\begin{align}
		\label{equ:rho^}
		\rho_K^{\alpha}=1+O(\varepsilon)
		,
	\end{align}     
	and 
	\begin{align}
		\label{equ:Wq=}
		\widetilde{W}_q(K)=\int_{\mathbb{S}^n}1 dx+O(\varepsilon)=|\mathbb{S}^n|+O(\varepsilon),
	\end{align}
	and 
	\begin{align}
		\label{equ:E(K)=}
		\widetilde{E}(K)=O(\varepsilon).
	\end{align}
\end{lemma}
\begin{proof}
	(1) Since $\rho_K\circ \mathscr{A}^*(x)=\sqrt{u_K^2(x)+|\nabla u_K|^2(x)}$,  $x\in \mathbb{S}^n$ (see \cite{Chen-Li-2021}*{eq. (2.4)} or \cite{Ding-Li-2025-JDE}*{Section 2.3}), we have 
	\begin{align*}
		\rho_K^2=&(1+w)^2+|\nabla w|^2=1+\left(2w+w^2+|\nabla w|^2\right)
		\\
		:=&1+g(w).
	\end{align*}
	Then, for $\varepsilon$ small enough so that $|g(w)|<1$, Taylor expansion gives
	\begin{align*}
		\rho_K^{\alpha}=&\left(1+g(w)\right)^\frac{\alpha}{2}
		\\
		=&1+\sum_{i=1}^{\infty}\binom{\frac{\alpha}{2}}{i}\cdot\left(g(w)\right)^{i}.
	\end{align*}
	Here $\binom{p}{i} := \frac{p(p-1)\cdots(p-i+1)}{i!}$ denotes the generalized binomial coefficient, for any real $p$ and integer $i \ge 0$. Since $|g(w)|\leq C_1\varepsilon$, we have $\left| \sum_{i=1}^{\infty}\binom{\frac{\alpha}{2}}{i}\cdot\left(g(w)\right)^{i}\right| \leq C_2\varepsilon$, with constant $C_1,C_2$. Then we obtain \eqref{equ:rho^}.
	\\
	(2)Since $K\in\mathcal{N}_{\varepsilon}$, we have $|w|,|w_i|,|w_{ij}|<O(\varepsilon)$ then 
	\begin{align}\label{equ:tau=1+O}
		\tau_{ij}[u_K]=w_{ij}+(1+w)\delta_{ij}=\delta_{ij}+O(\varepsilon),
	\end{align}
	and
	\begin{align}
		\label{sigma n=1+o}
		\sigma_n\left(\tau_{ij}[u_K]\right)=1+O(\varepsilon).
	\end{align}                                                                                                         
	Putting \eqref{equ:rho^} and \eqref{sigma n=1+o} into \eqref{equ:dual quermassintegral-C2},  we obtain \eqref{equ:Wq=}.
	\\
	(3) Taking $\varepsilon$ small enough so that $|g(w)|<1$, Taylor expansion yields 
	\begin{align}
		\label{equ:log rho}
		\log \rho_K=&\frac{1}{2}\log \left(1+g(w)\right)
		\nonumber
		\\
		=&\frac{1}{2}\left(
		\sum_{i=1}^{\infty}(-1)^{i+1}\cdot \frac{1}{i}\cdot \left(g(w)\right)^i
		\right).
	\end{align}
	Since $|g(w)|\leq C_3\varepsilon$, we have $\left|\frac{1}{2}\left(
	\sum_{i=1}^{\infty}(-1)^{i+1}\cdot \frac{1}{i}\cdot \left(g(w)\right)^i
	\right)\right| \leq C_4\varepsilon$, with constant $C_3,C_4$. Then we obtain \eqref{equ:E(K)=}.
\end{proof}


\begin{lemma}
	\label{lemma:int vEk<OE}
	Let $q\in\mathbb{R}$, and $K,L\in\mathcal{N}_{\varepsilon}$ be  $C^3$  convex bodies. If  $\widetilde{W}_q(K)=\widetilde{W}_q(L)$ for $q\neq 0$, then we have 
	\begin{align}
		\label{equ:int vEk<OE}
		\int_{\mathbb{S}^n}v\rho_K^{q-n-1}\sigma_n\left(\tau[u_t]\right) dx=-\frac{q-1}{2}\Vert v\Vert_{L^2}^2+\frac{1}{2}\Vert \nabla v\Vert_{L^2}^2+\mathscr{E}.
	\end{align}
	If $\widetilde{E}(K)=\widetilde{E}(L)$, then we have
	\begin{align}
		\label{equ:int v G(u) q=0}
		\int_{\mathbb{S}^n}v\rho_K^{-n-1}\sigma_n\left(\tau[u_t]\right) dx=\frac{                            1}{2}\Vert v\Vert_{L^2}^2+\frac{1}{2}\Vert \nabla v\Vert_{L^2}^2+\mathscr{E}.
	\end{align}
\end{lemma}
\begin{proof}
(1)	For any $t\in[0,1]$, we take a convex body $K_t=tL+(1-t)K$. Since $K,L\in\mathcal{N}_{\varepsilon}$, by Lemma \ref{lemma:C2-1-p-sum-C2-1}(i), we know $K_t\in\mathcal{N}_{\varepsilon}$ with support function $u_t:=tu_L+(1-t)u_K=u_K+tv$. Then we have $u_0=u_K$, $u_1=u_L$,  $\frac{d}{dt}u_t=v$ and $\frac{d}{dt}\tau[u_t]=\tau[v]$. We denote $r_t(x):=\rho_t\circ \mathscr{A}_{K_t}^*(x):=\sqrt{u_t^2(x)+|\nabla u_t|^2(x)}$, by \cite{Ding-Li-Trans}*{eq. (2.26)}, we have
	\begin{align}
		\label{equ:pf-lemma-alpha>-2}
		\frac{\partial \rho_t(\xi)}{\partial t}=\frac{\rho_t}{u_t}\frac{\partial u_t(x)}{\partial t}=\frac{\rho_t(\xi)}{u_t(x)}\cdot v(x),
	\end{align}	
	where $\xi=\mathscr{A}_{K_t}^*(x)$. 
	Then we can calculate that
	\begin{align}
		\frac{\partial }{\partial t}\widetilde{W}_q(K_t)\overset{\eqref{equ:dual quermassintegral}}{=}&\frac{\partial }{\partial t}\int_{\mathbb{S}^n}\rho_t^q(\xi) d\xi		
		=q\int_{\mathbb{S}^n}\rho_t^{q-1}(\xi) \frac{\partial\rho_t(\xi)}{\partial t} d\xi
		\nonumber
		\\
		\overset{\eqref{equ:pf-lemma-alpha>-2},\eqref{equ:dx=Adxi}}{=}		&
		q\int_{\mathbb{S}^n}r_t^{q-1} \cdot \frac{r_t}{u_t}\cdot v\cdot \frac{u_t\sigma_n(\tau[u_t])}{r_t^{n+1}} dx
		\nonumber
		\\
		=&q\int_{\mathbb{S}^n}v\cdot r_t^{q-n-1}\sigma_n\left(\tau[u_t]\right)dx\nonumber
		\\
		:=&q\int_{\mathbb{S}^n}vG(u_t)dx.
		\label{equ:pf-Lemma 3.2}
	\end{align}
	Here we denote \begin{align*}
		G(u_t)=r_t^{q-n-1}\sigma_n\left(\tau[u_t]\right).
	\end{align*}
	 Since 
	\begin{align}
		\label{equ:pt rt}
		\frac{\partial r_t}{\partial t}=&\frac{1}{2}\left(u_t^2+|\nabla u_t|^2\right)^{-\frac{1}{2}}\cdot\left(
		2u_t\cdot v+2\<\nabla u_t,\nabla v\>
		\right)
		\nonumber
		\\
		=&\frac{u_t v+\<\nabla u_t,\nabla v\>}{r_t},
	\end{align}
	and 
	\begin{align}
		\label{equ:pt sigma n}
		\frac{\partial \sigma_n(\tau[u_t])}{\partial t}=\sigma_n^{ij}(\tau[u_t])\cdot\left(v_{ij}+v\delta_{ij}\right).
	\end{align}
	Here $\sigma_n^{ij}(\tau[u_t])=\frac{\partial \sigma_n(\tau[u_t])}{\partial\tau_{ij}[u_t]}$. 
	Then we have
	\begin{align}
		\label{equ:pt G(u)}
		\frac{\partial}{\partial t}G(u_t)=&(q-n-1)r_t^{q-n-2}\cdot\frac{\partial r_t}{\partial t}\cdot \sigma_n\left(\tau[u_t]\right)+r_t^{q-n-1}\cdot\frac{\partial \sigma_n\left(\tau[u_t]\right)}{\partial t}\nonumber
		\\
		=&(q-n-1)r_t^{q-n-3}\cdot\left(u_t v+\<\nabla u_t,\nabla v\>\right)\cdot \sigma_n\left(\tau[u_t]\right)\nonumber
		\\&+r_t^{q-n-1}\cdot\sigma_n^{ij}(\tau[u_t])\cdot\left(v_{ij}+v\delta_{ij}\right).
	\end{align}
	 Using integration by parts, we can calculate that 
	 \begin{align}
	 	&\int_{\mathbb{S}^n}v\cdot  r_t^{q-n-1}\cdot\sigma_n^{ij}(\tau[u_t])\cdot v_{ij} dx
	 	\nonumber
	 	\\
	 	=&-\int_{\mathbb{S}^n}v_jv_i\cdot  r_t^{q-n-1}\cdot\sigma_n^{ij}(\tau[u_t]) +(q-n-1)vr_t^{q-n-2}(r_t)_i v_j\sigma_n^{ij}(\tau[u_t]) dx
	 	\nonumber
	 	\\
	 	=&-\int_{\mathbb{S}^n}v_jv_i\cdot  r_t^{q-n-1}\cdot\sigma_n^{ij}(\tau[u_t]) +(q-n-1)vr_t^{q-n-3}\<\nabla v,\nabla u_t\>\sigma_n(\tau[u_t]) dx,
	 	\label{equ:int vG-1}
	 \end{align}
	 where we use (see \cite{Huisken-Sinestrari-1999,LiHZ-2008})
	 \begin{align}
	 	\label{equ:Q_ij,j=0}
	 	&\sum_{j=1}^{n}\nabla_j\frac{\partial \sigma_n(\tau)}{\partial \tau_{ij}}=0,
	 	\\
	 	&\sum_{j=1}^{n}\frac{\partial \sigma_n(\tau)}{\partial \tau_{ij}} \tau_{jk}=\sigma_n(\tau)\delta_{ik}\label{equ:sigma n ij ik},
	 \end{align} 
	  and 
	 \begin{align}
	 	\label{equ:Di r}
	 	(r_t)_i=\nabla_i r_t=\frac{u_t \cdot  (u_t)_i+\sum_{k=1}^{n}(u_t)_{ki}(u_t)_k}{r_t}=\frac{\sum_{k=1}^{n}(u_t)_k\tau_{ki}[u_t]}{r_t}.
	 \end{align} 
	Using \eqref{equ:pf-Lemma 3.2}, \eqref{equ:pt G(u)} \eqref{equ:int vG-1}, we compute the second variation of $\widetilde{W}_q(K_t)$ as follows:
	\begin{align}\label{equ:ddW}
		\frac{\partial^2 \widetilde{W}_q(K_t)}{\partial t^2}=&q\int_{\mathbb{S}^n}v\left( (q-n-1)r_t^{q-n-3}\cdot\left(u_t v+\<\nabla u_t,\nabla v\>\right)\cdot \sigma_n\left(\tau[u_t]\right)\right.\nonumber
		\\&\left.+r_t^{q-n-1}\cdot\sigma_n^{ij}(\tau[u_t])\cdot\left(v_{ij}+v\delta_{ij}\right)\right)dx\nonumber
		\\
		=&q\int_{\mathbb{S}^n}(q-n-1)v^2 r_t^{q-n-3}u_t \cdot \sigma_n\left(\tau[u_t]\right)
		-v_jv_i\cdot  r_t^{q-n-1}\cdot\sigma_n^{ij}(\tau[u_t])\nonumber
		\\&+v^2r_t^{q-n-1}\cdot\sigma_n^{ij}(\tau[u_t])\cdot \delta_{ij}dx.
	\end{align}
	Since $K_t\in\mathcal{N}_{\varepsilon}$, from (see \cite{LiHZ-2008}*{eq. (2.27)})
	\begin{align*}
		\frac{\partial \sigma_{n}(\tau)}{\partial \tau_{ij}}=\frac{1}{(n-1)!} \sum_{i_1, \cdots, i_{n-1}; j_1, \cdots, j_{n-1}} \delta_{j_1 \cdots j_{n-1},j}^{i_1 \cdots i_{n-1},i}\tau_{i_1 j_1} \cdots \tau_{i_{n-1} j_{n-1}}
	\end{align*}
	 and \eqref{equ:tau=1+O} we can verify that
	\begin{align}\label{equ:sigmg ij=1+O}
	\sigma_n^{ij}(\tau[u_t])	=\delta_{ij}+O(\varepsilon).
	\end{align}
	Combining the above with \eqref{equ:rho^} and \eqref{sigma n=1+o}, and substituting them into \eqref{equ:ddW}, we obtain
	\begin{align}
		\label{equ:ddw=1+O}
		&\frac{\partial^2 \widetilde{W}_q(K_t)}{\partial t^2}\nonumber
		\\=&q\int_{\mathbb{S}^n}(q-n-1)v^2 \left(1+O(\varepsilon)\right)
		-\sum_{i,j=1}^{n}v_jv_i\left(\delta_{ij}+O(\varepsilon)\right)
		+v^2\left(1+O(\varepsilon)\right)
		\cdot ndx\nonumber
		\\
		=&q(q-1)\Vert v\Vert_{L^2}^2-q\Vert \nabla v\Vert_{L^2}^2+\mathscr{E}.
			\end{align}
	From \eqref{equ:pf-Lemma 3.2}, \eqref{equ:ddw=1+O},  Taylor expansion (with integral form of the remainder), and $\widetilde{W}_q(K)=\widetilde{W}_q(L)$, it follows that
	\begin{align*}
		0=&\widetilde{W}_q(K_1)-\widetilde{W}_q(K_0)=\left.\frac{\partial \widetilde{W}_q(K_t)}{\partial t}\right|_{t=0} +\int_0^1(1-t)\frac{\partial^2 \widetilde{W}_q(K_t)}{\partial t^2}dt\nonumber
		\\
		=&q\int_{\mathbb{S}^n}vG(u_0)dx+\left(q(q-1)\Vert v\Vert_{L^2}^2-q\Vert \nabla v\Vert_{L^2}^2+\mathscr{E}\right)\int_0^1(1-t)dt\nonumber
		\\
		=&q\int_{\mathbb{S}^n}v\rho_K^{q-n-1}\sigma_n\left(\tau[u_K]\right)dx+\frac{q(q-1)}{2}\Vert v\Vert_{L^2}^2-\frac{q}{2}\Vert \nabla v\Vert_{L^2}^2+\mathscr{E}
		.
	\end{align*}
This yields \eqref{equ:int vEk<OE}, because $q\neq 0$.
\\
(2) Similarly, we compute that
\begin{align}
	\frac{\partial }{\partial t}\widetilde{E}(K_t)\overset{\eqref{equ: dual entropy}}{=}&\frac{\partial }{\partial t}\int_{\mathbb{S}^n}\log \rho_t(\xi) d\xi	
	=                           \int_{\mathbb{S}^n}\rho_t^{ -1}(\xi) \frac{\partial\rho_t(\xi)}{\partial t} d\xi
	\nonumber
	\\
	\overset{\eqref{equ:pf-lemma-alpha>-2},\eqref{equ:dx=Adxi}}{=}		&
	    \int_{\mathbb{S}^n}r_t        ^{-1} \cdot \frac{r_t}{u_t}\cdot v\cdot \frac{u_t\sigma_n(\tau[u_t])}{r_t^{n+1}} dx
	\nonumber
	\\
	=&\int_{\mathbb{S}^n}v\cdot r_t^{-n-1}\sigma_n\left(\tau[u_t]\right)dx\nonumber
	\\
	:=&\int_{\mathbb{S}^n}vF(u_t)dx.
	\label{equ:pf-Lemma 3.2-0}
\end{align}
Here we denote \begin{align*}
	F(u_t)=r_t^{                   -n-1}\sigma_n\left(\tau[u_t]\right).
\end{align*}
By \eqref{equ:pt rt} and \eqref{equ:pt sigma n}, we have
\begin{align}
	\label{equ:pt G(u)-q=0}
	\frac{\partial}{\partial t}F(u_t)=&(    -n-1)r_t^{-n-2}\cdot\frac{\partial r_t}{\partial t}\cdot \sigma_n\left(\tau[u_t]\right)+r_t^{-n-1}\cdot\frac{\partial \sigma_n\left(\tau[u_t]\right)}{\partial t}\nonumber
	\\
	=&(-n-1)r_t^{-n-3}\cdot\left(u_t v+\<\nabla u_t,\nabla v\>\right)\cdot \sigma_n\left(\tau[u_t]\right)\nonumber
	\\&+r_t^{-n-1}\cdot\sigma_n^{ij}(\tau[u_t])\cdot\left(v_{ij}+v\delta_{ij}\right).
\end{align}
Using integration by parts, we can calculate that 
\begin{align}
	&\int_{\mathbb{S}^n}v\cdot  r_t^{-n-1}\cdot\sigma_n^{ij}(\tau[u_t])\cdot v_{ij} dx
	\nonumber
	\\
	=&-\int_{\mathbb{S}^n}v_jv_i\cdot  r_t^{-n-1}\cdot\sigma_n^{ij}(\tau[u_t]) +(-n-1)vr_t^{-n-2}(r_t)_i v_j\sigma_n^{ij}(\tau[u_t]) dx
	\nonumber
	\\
	=&-\int_{\mathbb{S}^n}v_jv_i\cdot  r_t^{-n-1}\cdot\sigma_n^{ij}(\tau[u_t]) +(-n-1)vr_t^{-n-3}\<\nabla v,\nabla u_t\>\sigma_n(\tau[u_t]) dx,
	\label{equ:int vG-1-q=0}
\end{align}
where we use \eqref{equ:Q_ij,j=0}, \eqref{equ:sigma n ij ik}, and \eqref{equ:Di r}.

Using \eqref{equ:pf-Lemma 3.2-0}, \eqref{equ:pt G(u)-q=0} \eqref{equ:int vG-1-q=0}, we compute the second variation of $\widetilde{E}(K_t)$ as follows:
\begin{align}\label{equ:ddW-q=0}
	\frac{\partial^2 \widetilde{E}(K_t)}{\partial t^2}=&\int_{\mathbb{S}^n}v\left( (-n-1)r_t^{-n-3}\cdot\left(u_t v+\<\nabla u_t,\nabla v\>\right)\cdot \sigma_n\left(\tau[u_t]\right)\right.\nonumber
	\\&\left.+r_t^{-n-1}\cdot\sigma_n^{ij}(\tau[u_t])\cdot\left(v_{ij}+v\delta_{ij}\right)\right)dx\nonumber
	\\
	=&\int_{\mathbb{S}^n}(-n-1)v^2 r_t^{-n-3}u_t \cdot \sigma_n\left(\tau[u_t]\right)
	-v_jv_i\cdot  r_t^{-n-1}\cdot\sigma_n^{ij}(\tau[u_t])\nonumber
	\\&+v^2r_t^{-n-1}\cdot\sigma_n^{ij}(\tau[u_t])\cdot \delta_{ij}dx.
\end{align}
Using \eqref{equ:sigmg ij=1+O}, \eqref{equ:rho^}, and \eqref{sigma n=1+o}, and substituting them into \eqref{equ:ddW-q=0}, we obtain
\begin{align}
	\label{equ:ddw=1+O-0}
	&\frac{\partial^2 \widetilde{E}(K_t)}{\partial t^2}\nonumber
	\\=&    \int_{\mathbb{S}^n}(-n-1)v^2 \left(1+O(\varepsilon)\right)
	-\sum_{i,j=1}^{n}v_jv_i\left(\delta_{ij}+O(\varepsilon)\right)
	+v^2\left(1+O(\varepsilon)\right)
	\cdot ndx\nonumber
	\\
	=&-\Vert v\Vert_{L^2}^2-\Vert \nabla v\Vert_{L^2}^2+\mathscr{E}.
\end{align}
From \eqref{equ:pf-Lemma 3.2-0}, \eqref{equ:ddw=1+O-0},  Taylor expansion (with integral form of the remainder), and $\widetilde{E}(K)=\widetilde{E}(L)$, it follows that
\begin{align*}
	0=&\widetilde{E}(K_1)-\widetilde{E}(K_0)=\left.\frac{\partial \widetilde{E}(K_t)}{\partial t}\right|_{t=0} +\int_0^1(1-t)\frac{\partial^2 \widetilde{E}(K_t)}{\partial t^2}dt\nonumber
	\\
	=&\int_{\mathbb{S}^n}vF(u_0)dx+\left(        -\Vert v\Vert_{L^2}^2-\Vert \nabla v\Vert_{L^2}^2+\mathscr{E}\right)\int_0^1(1-t)dt\nonumber
	\\
	=&\int_{\mathbb{S}^n}v\rho_K^{-n-1}\sigma_n\left(\tau[u_K]\right)dx-\frac{1}{2}\Vert v\Vert^2_{L^2}-\frac{1}{2}\Vert \nabla v\Vert^2_{L^2}+\mathscr{E}
	.
\end{align*}
That is \eqref{equ:int v G(u) q=0}.

	Then we complete the proof.
\end{proof}

Consequently, we can express the deficit of dual $L^p$-Minkowski inequality  in the form of the following lemmas:
\begin{lemma}
	\label{lemma:I_k-w-barw-O(e)}
	Let $q\in\mathbb{R}$, and let $K,L\in\mathcal{N}_{\varepsilon}$ be $C^3$ convex bodies. If  $\widetilde{W}_q(K)=\widetilde{W}_q(L)$ for $q\neq 0$ (and if $\widetilde{E}(K)=\widetilde{E}(L)$ when $q=0$), then for $p=0$, we have
	\begin{align}    
		\label{equ:deficit=p=0}
		\int_{\mathbb{S}^n}\left(\log(\frac{u_L}{u_K})\right)u_K\rho_K^{q-n-1}\sigma_n(\lambda_K)~dx
		=-\frac{q}{2}\Vert v\Vert_{L^2}^2
		+\frac{1}{2}\Vert\nabla v\Vert_{L^2}^2+\mathscr{E},
	\end{align}
	and for $p\neq 0$ we have
	\begin{align}
		\label{equ:deficit=p!=0}
		\frac{1}{p}\left(\int_{\mathbb{S}^n}
		u_L^pu_K^{1-p}\rho_K^{q-n-1}\sigma_n(\lambda_K)~dx-
		\int_{\mathbb{S}^n}
		u_K\rho_K^{q-n-1}\sigma_n(\lambda_K)~dx \right)\nonumber
		\\=\frac{p-q}{2}\Vert v\Vert_{L^2}^2
		+\frac{1}{2}\Vert\nabla v\Vert_{L^2}^2+\mathscr{E}.
	\end{align}
\end{lemma}
\begin{remark}
For $q\neq 0$, the second term inside the parentheses on the left-hand side of \eqref{equ:deficit=p!=0} is $\widetilde{W}_q(K)$, that is 
\begin{align*}
 \frac{1}{p}\left(\int_{\mathbb{S}^n}
 u_L^pu_K^{1-p}\rho_K^{q-n-1}\sigma_n(\lambda_K)~dx-\widetilde{W}_q(K) \right)\nonumber
 =\frac{p-q}{2}\Vert v\Vert_{L^2}^2
 +\frac{1}{2}\Vert\nabla v\Vert_{L^2}^2+\mathscr{E}.
\end{align*}
For $q=0$, by \eqref{equ:dx=Adxi}, the second term in the parentheses on the left-hand side of \eqref{equ:deficit=p!=0} equals $|\mathbb{S}^n|$, that is,
\begin{align*}
	\frac{1}{p}\left(\int_{\mathbb{S}^n}
	u_L^pu_K^{1-p}\rho_K^{-n-1}\sigma_n(\lambda_K)~dx-|\mathbb{S}^n|\right)\nonumber
	=\frac{p}{2}\Vert v\Vert_{L^2}^2
	+\frac{1}{2}\Vert\nabla v\Vert_{L^2}^2+\mathscr{E}.
\end{align*}
\end{remark}
	\begin{proof}[\textbf{Proof of Lemma  \ref{lemma:I_k-w-barw-O(e)}}]
(1)	We take $\varepsilon$ small enough so that $\tau[u_K]>0,\tau[u_L]>0$, $\frac{1}{2}<u_K<\frac{3}{2}$, and $\left|\frac{v}{u_K}\right|<1$.
	By  Taylor expansion, we have
	\begin{align*}
		\log\left(\frac{u_L}{u_K}\right)=&\log\left(1+\frac{v}{u_K}\right)
		\\=&\frac{v}{u_K}-\frac{1}{2}\left(\frac{v}{u_K}\right)^2+\sum_{i=3}^{\infty}(-1)^{i+1}\frac{1}{i}\left(\frac{v}{u_K}\right)^i
		\\=&\frac{v}{u_K}-\frac{1}{2}\left(\frac{v}{u_K}\right)^2+O(\varepsilon)\left(\frac{v}{u_K}\right)^2.
	\end{align*}
	The last equality above follows from the estimate $$\left|\sum_{i=3}^{\infty}(-1)^{i+1}\frac{1}{i}\left(\frac{v}{u_K}\right)^i\right|\leq O(\varepsilon)\left(\frac{v}{u_K}\right)^2.$$
Since $\widetilde{W}_q(K)=\widetilde{W}_q(L)$ (or $\widetilde{E}(K)=\widetilde{E}(L)$ for $q=0$), combining the above with Lemma \ref{lemma:int vEk<OE}, \eqref{equ:rho^}, and \eqref{sigma n=1+o}, we obtain
	\begin{align*}       &\int_{\mathbb{S}^n}\left(\log(\frac{u_L}{u_K})\right)u_K\rho_K^{q-n-1}\sigma_n(\lambda_K)~dx
		\\
		=&\int_{\mathbb{S}^n}
		v\rho_K^{q-n-1}\sigma_n(\lambda_K)dx-\frac{1}{2}\int_{\mathbb{S}^n}
		v^2u_K^{-1}\rho_K^{q-n-1}\sigma_n(\lambda_K)dx+O(\varepsilon)\int_{\mathbb{S}^n}
		v^2u_K^{-1}\rho_K^{q-n-1}\sigma_n(\lambda_K)dx
		\\
		=&-\frac{q-1}{2}\Vert v\Vert^2_{L^2}+\frac{1}{2}\Vert \nabla v\Vert^2_{L^2}+\mathscr{E}
		-\frac{1}{2}\Vert v\Vert_{L^2}^2\left(1+O(\varepsilon)\right)+O(\varepsilon)\Vert v\Vert_{L^2}^2
		\\
		=&
		-\frac{q}{2}\Vert v\Vert_{L^2}^2
		+\frac{1}{2}\Vert\nabla v\Vert_{L^2}^2+O(\varepsilon)\Vert  v\Vert_{L^2}^2+O(\varepsilon)\Vert \nabla v\Vert_{L^2}^2,
	\end{align*}
	where we use the fact that $u_K^{-1}=\frac{1}{1+O(\varepsilon)}=1+O(\varepsilon)$. Then we obtain \eqref{equ:deficit=p=0}.
	\\
	(2)
	When $p\neq 0$, from  Taylor expansion, we have
	\begin{align*}
		\left(\frac{u_L}{u_K}\right)^p=&\left(1+\left(\frac{v}{u_K}\right)\right)^p
		\\
		=&1+p\left(\frac{v}{u_K}\right)+\frac{p(p-1)}{2}\left(\frac{v}{u_K}\right)^2
		+\sum_{i=3}^{\infty}\binom{p}{i}\left(\frac{v}{u_K}\right)^i
		\\
		=&1+p\left(\frac{v}{u_K}\right)+\frac{p(p-1)}{2}\left(\frac{v}{u_K}\right)^2+O(\varepsilon)\left(\frac{v}{u_K}\right)^2,
	\end{align*}
	where the last equality above follows from the estimate $$\left|\sum_{i=3}^{\infty}\binom{p}{i}\left(\frac{v}{u_K}\right)^i\right|\leq O(\varepsilon)\left(\frac{v}{u_K}\right)^2.$$
	Together with Lemma \ref{lemma:int vEk<OE}, we obtain
	\begin{align*}
		&\frac{1}{p}\left(\int_{\mathbb{S}^n}
		u_L^pu_K^{1-p}\rho_K^{q-n-1}\sigma_n(\lambda_K)~dx-\int_{\mathbb{S}^n}
		u_K\rho_K^{q-n-1}\sigma_n(\lambda_K)~dx \right)
		\\
		=&\frac{1}{p}\int_{\mathbb{S}^n}
		u_K\left(\left(1+\left(\frac{v}{u_K}\right)\right)^p-1\right)\rho_K^{q-n-1}\sigma_n(\lambda_K)
		~dx
		\\
		=&\frac{1}{p}\int_{\mathbb{S}^n}
		u_K\left(p\left(\frac{v}{u_K}\right)+\frac{p(p-1)}{2}\left(\frac{v}{u_K}\right)^2+O(\varepsilon)\left(\frac{v}{u_K}\right)^2\right)\rho_K^{q-n-1}\sigma_n(\lambda_K)
		~dx
		\\
		=&
		-\frac{q-1}{2}\Vert v\Vert_{L^2}^2
		+\frac{1}{2}\Vert\nabla v\Vert_{L^2}^2+\mathscr{E}+\frac{p-1}{2}\Vert v\Vert_{L^2}^2+O(\varepsilon)\Vert  v\Vert_{L^2}^2.
	\end{align*}
	Then we obtain \eqref{equ:deficit=p!=0}.
\end{proof}

Motivated by \cite[eq. (73)]{Cicalese-Leonardi-2012} and \cite[Lemma 3.4]{Ding-Gao-Li-BM1}, we establish the following lower bound for  $\Vert\nabla(w-\bar{w})\Vert_{L^2}^2$. 

An analogous estimate was established in \cite[Lemma 3.4]{Ding-Gao-Li-BM1} under the assumption that the two convex bodies have equal intrinsic volumes. In the present setting, where the two bodies have equal $\widetilde{W}_q$ or $\widetilde{E}$, the only difference lies in the coefficient $a_0$. More precisely, this estimate now follows from the equality of $W_q$ or $E$, as shown in \eqref{equ:a_0=OE}. The remainder of the argument is identical to that in \cite[Lemma 3.4]{Ding-Gao-Li-BM1}. Nevertheless, for the sake of completeness, we provide the proof below.
	\begin{lemma}
	\label{lemma:DU>(2n+1)u}
	Let $p,q\in \mathbb{R}$, and $K,L\in \mathcal{N}_{\varepsilon}$ be    convex bodies  with support function $u_K=1+w,u_L=1+\bar{w}\in C^3(\mathbb{S}^n)$, respectively. Assume  $\widetilde{W}_q(K)=\widetilde{W}_q(L)$ with $q\neq 0$ (or $\widetilde{E}(K)=\widetilde{E}(L)$ when $q=0$), and the $p$-weighted  Steiner points satisfy  $\mathcal{S}_p(K)=\mathcal{S}_p(L) $.
	Then we have
	\begin{align}
		\label{equ:equ:DU>(2n+1)u}
		\Vert\nabla(u_K-u_L)\Vert_{L^2}^2\geq2(n+1)\Vert u_K-u_L\Vert_{L^2}^2+\mathscr{E}.
	\end{align}
\end{lemma}
\begin{proof}
	Denote by $\{Y_{\ell}\}$ the spherical harmonics which form an orthonormal basis for $L^2(\mathbb{S}^n)$, i.e., 
	\begin{align}\label{equ:yiyj}
		\int_{\mathbb{S}^n}Y_{k}Y_{\ell} dx=\delta_{k \ell}, \quad \text{for\ } k,\ell=0,1,2,\cdots,
	\end{align}
	and
	\begin{align}
		\label{equ:Laplave Y=lambda Y}
		-\Delta_{\mathbb{S}^n} Y_{\ell}	=\lambda_{\ell}Y_{\ell}, \quad \text{for\ }\ell =0,1,2,\cdots,
	\end{align}
	where $\lambda_{\ell}=\ell(\ell+n-1)$ are the eigenvalues of the negative Laplace-Beltrami operator $-\Delta_{\mathbb{S}^n}$ on the sphere, with corresponding eigenfunctions being the spherical harmonics of degree $\ell$. 
	We write
	\begin{align}
		\label{equ:w-w-}
		w-\bar{w}=\sum_{\ell=0}^{\infty}a_{\ell}Y_{\ell}.
	\end{align}
	Recall that
	\begin{align*}
		Y_0=\frac{1}{\sqrt{|\mathbb{S}^n|}},\quad Y_1(x)=\<x,e\>,
	\end{align*}
	for a suitably chosen constant vector $e$. Here $|\mathbb{S}^n|$ denotes the surface area of $\mathbb{S}^n$.
	
	We have
	\begin{align*}
		a_0=\int_{\mathbb{S}^n}( w-\bar{w})Y_0~dx=\frac{1}{\sqrt{|\mathbb{S}^n|}}\int_{\mathbb{S}^n}( w-\bar{w})~dx,
	\end{align*}
	and
	\begin{align*}
		a_1=\int_{\mathbb{S}^n}( w-\bar{w})Y_1~dx=\left\langle \mathcal{S}(K)-\mathcal{S}(L),e\right\rangle. 
	\end{align*}
	We write $u_t=u_K+tv$, $0\leq t\leq 1$, and 
	\begin{equation*}
		r_t=(u_t^2+|\nabla u_t|^2)^{1/2},\qquad
		G_t=r_t^{q-n-1}\det\tau[u_t],\qquad
		\overline G(x)=\int_0^1G_t(x)\,dt.
	\end{equation*}
	Since 	 either $\widetilde{W}_q(K)=\widetilde{W}_q(L)$ for $q\neq 0$ or $\widetilde{E}(K)=\widetilde{E}(L)$, from \eqref{equ:pf-Lemma 3.2} and \eqref{equ:pf-Lemma 3.2-0}, we have
	$$0=\widetilde{W}_q(L)-\widetilde{W}_q(K)=\int _0^1q\int_{\mathbb S^n}vG_tdx, \text{\ for\ } q\ne 0,$$ 
	and  $$0=\widetilde{E}(L)-\widetilde{E}(K)=\int_{\mathbb S^n}vG_tdx, \text{\ for\ }q=0.$$
	That is
	$\int_{\mathbb S^n}v\overline G\,dx=0$ in both cases.
	Because $u_t$ is uniformly $C^2$-close to $1$,
	$\|G_t-1\|_{L^\infty}\leq C_{n,q}\varepsilon$ uniformly in $t$.
	Then Cauchy--Schwarz inequality gives
	\begin{equation*}
		|a_0|=\frac{1}{\sqrt{ |\mathbb{S}^n|}}\left|\int_{\mathbb S^n}v(1-\overline G)\,dx\right|
		\leq C_{n,q}\varepsilon\|v\|_{L^2}.
	\end{equation*}
	Then
	\begin{align}\label{equ:a_0=OE}
		a_0^2=\mathscr{E}.
	\end{align}
	For any $p\in\mathbb R$, set $\Phi_p(z)=(z^p-1)/p$ for $p\ne0$
	and $\Phi_0(z)=\log z$. Since $\Phi_p'(z)=z^{p-1}$ for $z>0$, we have
	\begin{equation*}
		\Phi_p(u_L)-\Phi_p(u_K)=vB,\qquad
		B(x)=\int_0^1 u_t(x)^{p-1}\,dt.
	\end{equation*}
	By the definition of the $p$-weighted Steiner point, $\mathcal{S}_p(K)=\mathcal{S}_p(L)$ gives $\int_{\mathbb S^n}(\Phi_p(u_L)-\Phi_p(u_K))x\,dx=0$.
	Thus, also in the logarithmic case,
	\begin{equation*}
		\int_{\mathbb S^n}vB Y_1\,dx=0.
	\end{equation*}
	As $u_t=1+O(\varepsilon)$ uniformly, $\|1-B\|_{L^\infty}\leq
	C_p\varepsilon$. Then we have
	\begin{equation*}
		|a_1|=\left| \int_{\mathbb{S}^n}vY_1dx\right| =\left| \int_{\mathbb{S}^n}v(1-B)Y_1dx\right|
		\leq C_p\varepsilon\|v\|_{L^2}.
	\end{equation*}
For any $p\in\mathbb{R}$,  we have
	\begin{align}\label{equ:a_1<OE}
		a_1^2\leq O(\varepsilon)\Vert w-\bar{w}\Vert_{L^2}^2.
	\end{align}
	
	By \eqref{equ:yiyj} and  \eqref{equ:w-w-} we have  $$\Vert w-\bar{w}\Vert_{L^2}^2=\sum_{\ell=0}^{\infty}a_{\ell}^2.$$  
	From \eqref{equ:yiyj}$\thicksim$\eqref{equ:w-w-} and  integration by parts, we can calculate that
	\begin{align*}
		\Vert\nabla( w-\bar{w})\Vert_{L^2}^2=&-\int_{\mathbb{S}^n}(w-\bar{w})\cdot \Delta_{\mathbb{S}^n}(w-\bar{w}) dx
		\\
		\overset{\eqref{equ:w-w-}}{=}&-\int_{\mathbb{S}^n}
		\left(\sum_{\ell=0}^{\infty}a_{\ell}Y_{\ell}\right)\cdot
		\left( \sum_{k=0}^{\infty}a_{k}\Delta_{\mathbb{S}^n}Y_{k}\right) dx
		\\
		\overset{\eqref{equ:Laplave Y=lambda Y}}{=}&\int_{\mathbb{S}^n}
		\left(\sum_{\ell=0}^{\infty}a_{\ell}Y_{\ell}\right)\cdot
		\left( \sum_{k=0}^{\infty}a_k\lambda_{k}Y_{k}\right) dx
		\\
		\overset{\eqref{equ:yiyj}}{=}& \sum_{\ell=0}^{\infty}\lambda_{\ell}a_{\ell}^2.
	\end{align*}
	Noting that $\lambda_{\ell}=\ell(\ell+n-1)\geq 2(n+1)$ when $\ell\geq 2$, we derive
	\begin{align*}
		\Vert\nabla( w-\bar{w})\Vert_{L^2}^2=&\sum_{\ell=0}^{\infty}\lambda_{\ell}a_{\ell}^2    =\sum_{\ell=2}^{\infty}\lambda_{\ell}a_{\ell}^2+na_1^2
		\\
		\geq& 2(n+1)\sum_{\ell=2}^{\infty}a_{\ell}^2+na_1^2
		\\
		=&2(n+1)\sum_{\ell=0}^{\infty}a_{\ell}^2-2(n+1)a_0^2-(n+2)a_1^2
		\\
		\geq&2(n+1)\Vert w-\bar{w}\Vert_{L^2}^2+\mathscr{E}.
	\end{align*}
	In the last equality, we used \eqref{equ:a_0=OE} and \eqref{equ:a_1<OE}.
\end{proof}

\section{Stability of dual $L^p$-Minkowski inequality}\label{sec 4}
Based on the computations in Section \ref{sec 3}, this section aims to establish the stability of dual  $L^p$-Minkowski inequality in terms of different distances. We prove Theorems \ref{thm:p-M-L2-Vk=Vk}$\thicksim$\ref{thm:p-M-L2-F-Vk=Vk}, 
as well as Corollaries \ref{thm:dual p-M-L2} and \ref{thm:dual p-M-al}.
\subsection{In terms of $L^2$ distance}
The following theorem, which is equivalent to Corollary \ref{thm:dual p-M-L2} in the case $q \neq 0$, is proved first.
\begin{theorem}
	\label{thm:p-M-L2-Vk=Vk}
	Let $K, L \in \mathcal{K}^{n+1}$ be $C^3$ convex bodies and $p,q$ be real numbers with  $p > q - 2n - 2, q\neq 0$. 
	For all $\eta>0$, there exists $\varepsilon>0$  such that  if 
	$K,L\in\mathcal{N}_{\varepsilon}$, $\widetilde{W}_q(K)=\widetilde{W}_q(L)$, and $K,L$ satisfy $S_p(K)=S_p(L)={o}$, then for $p=0$, we have
	\begin{align}\label{equ:p-M-Vk=p=0}
		\int_{\mathbb{S}^n}\left(\log(\frac{u_L}{u_K})\right)u_K\rho_K^{q-n-1}\sigma_n(\lambda_K) ~dx
		\geq \left(\frac{2n+2-q}{2}-\eta\right)\delta_2^2(K,L),
	\end{align}
	and for $p\neq 0$, we have
	\begin{align}\label{equ:p-M-Vk=p!=0}
		\frac{1}{p}\left(\int_{\mathbb{S}^n}u_L^pu_K^{1-p}\rho_K^{q-n-1}\sigma_n(\lambda_K)~dx-\widetilde{W}_q(K)\right)
		\geq \left(\frac{2n+2+p-q}{2}-{\eta}\right)\delta_2^2(K,L).
	\end{align}
\end{theorem} 

\begin{proof}[\textbf{Proof of Theorem \ref{thm:p-M-L2-Vk=Vk}}]
	(1).	Since 
	$$-\frac{\Vert \nabla w\Vert_{L^2}^2}{2}-\frac{\Vert \nabla\bar{w}\Vert_{L^2}^2}{2}\leq\int_{\mathbb{S}^n} \sum_{i=1}^n w_i\bar{w}_i ~dx\leq \frac{\Vert \nabla w\Vert_{L^2}^2}{2}+\frac{\Vert \nabla\bar{w}\Vert_{L^2}^2}{2},$$
	and
	$$-\frac{\Vert w\Vert_{L^2}^2}{2}-\frac{\Vert  \bar{w}\Vert_{L^2}^2}{2}\leq\int_{\mathbb{S}^n}  w\bar{w}~dx\leq \frac{\Vert  w\Vert_{L^2}^2}{2}+\frac{\Vert \bar{w}\Vert_{L^2}^2}{2},$$
	$\mathscr{E}$ can also be expressed as:
	\begin{align*}
		\mathscr{E} =O(\varepsilon)\Vert w-\bar{w}\Vert_{L^2}^2+O(\varepsilon)\Vert \nabla w-\nabla\bar{w}\Vert_{L^2}^2.
	\end{align*}
	Substituting the above expression  and \eqref{equ:equ:DU>(2n+1)u} into \eqref{equ:deficit=p=0} yields  
	\begin{align}   
		&\int_{\mathbb{S}^n}\left(\log(\frac{u_L}{u_K})\right)u_K\rho_K^{q-n-1}\sigma_n(\lambda_K)~dx		\nonumber
		\\
		=&\left(-\frac{q}{2}+O(\varepsilon)\right)\Vert\bar{w}-w\Vert_{L^2}^2
		+\left(\frac{1}{2}+O(\varepsilon)\right)\Vert\nabla\bar{w}-\nabla w\Vert_{L^2}^2\nonumber
		\\
		\geq& \left(-\frac{q}{2}+O(\varepsilon)\right)\Vert\bar{w}-w\Vert_{L^2}^2
		+2(n+1)\left(\frac{1}{2}+O(\varepsilon)\right)\Vert\bar{w}- w\Vert_{L^2}^2\nonumber
		\\
		=&\left(\frac{2n+2-q}{2}+O(\varepsilon)\right)\delta_2^2(K,L).\label{equ:pf-thm:p-M-L2-Vk=Vk-1}
	\end{align}
	\\	
	(2).  For  $p\neq 0$, the analysis proceeds analogously to the case $p=0$. Putting \eqref{equ:equ:DU>(2n+1)u} into \eqref{equ:deficit=p!=0} yields  
	\begin{align}
		&\frac{1}{p}\left(\int_{\mathbb{S}^n}u_L^pu_K^{1-p}\rho_K^{q-n-1}\sigma_n(\lambda_K)~dx-\widetilde{W}_q(K)\right)\nonumber
		\\
		=&\left(\frac{p-q}{2}+O(\varepsilon)\right)\Vert\bar{w}-w\Vert_{L^2}^2
		+\frac{2(n+1)}{2}\left(\frac{1}{2(n+1)}+O(\varepsilon)\right)\Vert\nabla\bar{w}-\nabla w\Vert_{L^2}^2\nonumber
		\\
		\geq&\left(\frac{p-q}{2}+O(\varepsilon)\right)\Vert\bar{w}-w\Vert_{L^2}^2
		+(n+1)\left(1+O(\varepsilon)\right)\Vert \bar{w}-  w\Vert_{L^2}^2\nonumber
		\\
		=& \left(\frac{2n+2+p-q}{2}+O(\varepsilon)\right)\delta_2^2(K,L).
		\label{equ:pf-thm:p-M-L2-Vk=Vk-2}
	\end{align}
\end{proof}
Now, Corollary \ref{thm:dual p-M-L2} for $q \neq 0$ follows directly from Theorem \ref{thm:p-M-L2-Vk=Vk}, while the case $q = 0$ follows from Lemma \ref{lemma:I_k-w-barw-O(e)}  and Lemma \ref{lemma:DU>(2n+1)u}.
\begin{corollary}
	\label{thm:dual p-M-L2} 
	Suppose $K, L \in \mathcal{K}^{n+1}$ are $C^3$ convex bodies, $p$ and $q$ are real numbers with  $p > q - 2n - 2$.
	For all $\eta>0$, there exists $\varepsilon>0$  such that if 
	$K,L\in\mathcal{N}_{\varepsilon}$, and $K,L$  satisfy  $S_p(K)=S_p(L)={o}$, then for $p\neq 0, q\neq 0$, we have
	\begin{align}\label{equ:dual p-M-L2}
		&\frac{1}{p}\left(\int_{\mathbb{S}^n}u_L^pu_K^{1-p}\rho_K^{q-n-1}\sigma_n(\lambda_K)~dx-\widetilde{W}_q(K)^{1-\frac{p}{q}}\widetilde{W}_q(L)^{\frac{p}{q}}\right)\nonumber
		\\
		\geq& N^{q-p-2}\left(\frac{2n+2+p-q}{2}-{\eta}\right)\delta_2^2( K,NL),
	\end{align}
	and for $p=0, q\neq 0$, we have
	\begin{align}\label{equ:dual p-M-L2-p=0}
		&\int_{\mathbb{S}^n}\left(\log(\frac{u_L}{u_K})\right)u_K\rho_K^{q-n-1}\sigma_n(\lambda_K) ~dx
		-\frac{1}{q}\widetilde{W}_q(K)\log\left(\frac{\widetilde{W}_q(L)}{\widetilde{W}_q(K)}\right)\nonumber
		\\
		\geq& N^{q-2}\left(\frac{2n+2-q}{2}-\eta\right)\delta_2^2( K,NL),
	\end{align}
	and for $p\neq 0, q=0$, 
	we have
	\begin{align}\label{equ:dual p-M-L2-q=0}
		&\frac{1}{p}\left(\int_{\mathbb{S}^n}u_L^pu_K^{1-p}\rho_K^{-n-1}\sigma_n(\lambda_K)~dx-|\mathbb{S}^n|e^{\frac{p\left(\widetilde{E}(L)-\widetilde{E}(K)\right)}{|\mathbb{S}^n|}}\right)
		\nonumber
		\\
		\geq& N_0^{-p-2}\left(\frac{2n+2+p}{2}-{\eta}\right)\delta_2^2( K,N_0L),
	\end{align}
	and for $p=0,q=0$, we have 
	\begin{align}
		\label{equ:dual p-M-L2-p=0-q=0}
		&\int_{\mathbb{S}^n}\log\left(\frac{u_L}{u_K}\right)u_K\rho_K^{-n-1}\sigma_n(\lambda_K)~dx+\widetilde{E}(K)-\widetilde{E}(L)\nonumber
		\\
		\geq &N_0^{-2}\left((n+1)-\eta\right)\delta_2^2(K,N_0L).
	\end{align}
	Here 
	$N=\left(\frac{\widetilde{W}_q(K)}{\widetilde{W}_q(L)}\right)^{\frac{1}{q}}$ and $N_0=e^{\frac{\widetilde{E}(K)-\widetilde{E}(L)}{|\mathbb{S}^n|}}$.
\end{corollary} 
\begin{remark}
	(1) When $p>0=q$, since $e^x\geq x+1$, inequality \eqref{equ:dual p-M-L2-q=0} reduces to
	\begin{align*}
		&\frac{1}{p}\left(\int_{\mathbb{S}^n}u_L^pu_K^{1-p}\rho_K^{-n-1}\sigma_n(\lambda_K)~dx-|\mathbb{S}^n|\right)+\widetilde{E}(K)-\widetilde{E}(L)
		\nonumber
		\\
		\geq
		&\frac{1}{p}\left(\int_{\mathbb{S}^n}u_L^pu_K^{1-p}\rho_K^{-n-1}\sigma_n(\lambda_K)~dx-|\mathbb{S}^n|e^{\frac{p\left(\widetilde{E}(L)-\widetilde{E}(K)\right)}{|\mathbb{S}^n|}}\right)
		\nonumber
		\\
		\geq& N_0^{-p-2}\left(\frac{2n+2+p}{2}-{\eta}\right)\delta_2^2( K,N_0L).
	\end{align*}
	which is the stability of inequality in \cite{Xi-Zhang-2022}*{eq. (4.4)}. Similarly, \eqref{equ:dual p-M-al-q=0} in Corollary \ref{thm:dual p-M-al} is also the stability  of the same inequality.

	(2) $\delta_2^2(K,NL) \geq 0$ implies that the dual $L^p$-Minkowski inequality also holds for $p > q-2n-2$ (see Theorem \ref{thm:p-M-inequ}), and equality holds if and only if $K=NL$.
\end{remark}
\begin{proof}[\textbf{Proof of Corollary \ref{thm:dual p-M-L2}}] 	
	(1). We begin by considering the case $q \neq 0$. For $K$ and $L$, there exists an constant $N>0$, such that $\widetilde{K}:=\frac{1}{N}K$ satisfies $\widetilde{W}_q(\widetilde{K})=\widetilde{W}_q(L)$. Then we have
	$$\left(\frac{1}{N}\right)^{q}\widetilde{W}_q(K)=\widetilde{W}_q(L),$$
	which implies 
	\begin{align*}
		N=\left(\frac{\widetilde{W}_q(K)}{\widetilde{W}_q(L)}\right)^{\frac{1}{q}}.
	\end{align*}
	From \eqref{equ:Wq=}, we have $\widetilde{W}_q(K)= |\mathbb{S}^{n}|+O(\varepsilon)$, and $\widetilde{W}_q(L)=|\mathbb{S}^{n}|+O(\varepsilon)$, then $N=1+O(\varepsilon)$. Denote $u_{\widetilde{K}}=\frac{u_K}{N}:=1+\tilde{w}$ is the support function of $\widetilde{K}$, then we have $\Vert \tilde{w}\Vert_{C^2}=O(\varepsilon)$.
	One can check that  pair $(\widetilde{K},L)$ satisfies $S_p(\widetilde{K})=S_p(L)={o}$.

	By $K=N\widetilde{K}$ and Theorem \ref{thm:p-M-L2-Vk=Vk}, we have 
	\begin{align*}
		&\int_{\mathbb{S}^n}\left(\log(\frac{u_L}{u_K})\right)u_K\rho_K^{q-n-1}\sigma_n(\lambda_K) ~dx
		-\frac{1}{q}\widetilde{W}_q(K)\log\frac{\widetilde{W}_q(L)}{\widetilde{W}_q(K)}
		\\
		=&\int_{\mathbb{S}^n}\left(\log(\frac{u_L}{Nu_{\widetilde{K}}})\right)Nu_{\widetilde{K}}N^{q-1}\rho_{\widetilde{K}}^{q-n-1}\sigma_n(\lambda_{\widetilde{K}}) ~dx-\frac{1}{q}N^{q}\widetilde{W}_q(\widetilde{K})\log\frac{\widetilde{W}_q(L)}{N^{q}\widetilde{W}_q(\widetilde{K})}
		\\
		=&N^{q}\int_{\mathbb{S}^n}\left(\log(\frac{u_L}{u_{\widetilde{K}}})\right)u_{\widetilde{K}}\rho_{\widetilde{K}}^{q-n-1}\sigma_n(\lambda_{\widetilde{K}}) ~dx
		\\
		\geq& N^{q}\left(\frac{2n+2-q}{2}-\eta\right)\delta_2^2(\widetilde{K},L)
		.
	\end{align*}
	Since 
	\begin{align*}
		\delta_2^2(\widetilde{K},L)=&\int_{\mathbb{S}^n}(u_{\widetilde{K}}-u_L)^2~dx
		=\int_{\mathbb{S}^n}\left(\frac{1}{N}u_K-u_L\right)^2~dx
		=N^{-2}\int_{\mathbb{S}^n}(u_K-Nu_L)^2~dx
		\\
		=&N^{-2}\delta_2^2({K},NL).
	\end{align*}
	Then we obtain the inequality \eqref{equ:dual p-M-L2-p=0}.
	
	In the case $p \neq 0$, a similar argument using $K = N\widetilde{K}$ together with \eqref{equ:p-M-Vk=p!=0} from Theorem \ref{thm:p-M-L2-Vk=Vk} gives
	\begin{align*}
		&\frac{1}{p}\left(\int_{\mathbb{S}^n}u_L^pu_K^{1-p}\rho_K^{q-n-1}\sigma_n(\lambda_K)~dx-\widetilde{W}_q(K)^{\frac{q-p}{q}}\widetilde{W}_q(L)^{\frac{p}{q}}\right)
		\\
		=&\frac{1}{p}\left(\int_{\mathbb{S}^n}u_L^pu_{(N\widetilde{K})}^{1-p}\rho_{(N\widetilde{K})}^{q-n-1}\sigma_n(\lambda_{(N\widetilde{K})})~dx-\widetilde{W}_q({N\widetilde{K}})^{\frac{q-p}{q}}\widetilde{W}_q(L)^{\frac{p}{q}}\right)
		\\
		=&N^{q-p}\frac{1}{p}\left(\int_{\mathbb{S}^n}u_L^pu_{\widetilde{K}}^{1-p}\rho_{\widetilde{K}}^{q-n-1}\sigma_n(\lambda_{\widetilde{K}})~dx-\widetilde{W}_q({\widetilde{K}})^{\frac{q-p}{q}}\widetilde{W}_q(L)^{\frac{p}{q}}\right)
		\\
		\geq& N^{q-p-2}\left(\frac{2n+2+p-q}{2}-{\eta}\right)\delta_2^2( K,NL).
	\end{align*}    
	Then we obtain the inequality \eqref{equ:dual p-M-L2}.
	\\
	(2). We next turn to the case $q = 0$. We choose $N_0>0$ such that $\tilde{K} = \frac{1}{N_0}K$ satisfies $\widetilde{E}(\widetilde{K}) = \widetilde{E}(L)$, that is, 
	$$\log N_0 = \frac{\widetilde{E}(K) - \widetilde{E}(L)}{|\mathbb{S}^n|}.$$  From  \eqref{equ:E(K)=}, we have $\widetilde{E}(K)= O(\varepsilon)$, and $\widetilde{E}(L)=O(\varepsilon)$, then $N_0=1+O(\varepsilon)$. Denote $u_{\widetilde{K}}=\frac{u_K}{N_0}:=1+\tilde{w}$ is the support function of $\widetilde{K}$, then we have $\Vert \tilde{w}\Vert_{C^2}=O(\varepsilon)$. We can check that  pair $(\widetilde{K},L)$ satisfies $S_p(\widetilde{K})=S_p(L)={o}$.
	
	If $p \neq 0$, substituting \eqref{equ:equ:DU>(2n+1)u} into \eqref{equ:deficit=p!=0} with $q=0$ yields
	\begin{align*}
		&\frac{1}{p}\left(\int_{\mathbb{S}^n}
		u_L^pu_{\widetilde{K}}^{1-p}\rho_{\widetilde{K}}^{-n-1}\sigma_n(\lambda_{\widetilde{K}})~dx-|\mathbb{S}^n|\right)
		\nonumber
		\\
		=&\frac{p}{2}\Vert\bar{w}-\tilde{w}\Vert_{L^2}^2
		+\frac{1}{2}\Vert\nabla\bar{w}-\nabla \tilde{w}\Vert_{L^2}^2+\mathscr{E}
		\\
		\geq &\left(\frac{2n+2+p}{2}+O(\varepsilon)\right)\Vert\bar{w}-\tilde{w}\Vert_{L^2}^2.
	\end{align*}
	From the above inequality we obtain
	\begin{align*}
		&\frac{N_0^p}{p}\left(\int_{\mathbb{S}^n}u_L^pu_K^{1-p}\rho_K^{-n-1}\sigma_n(\lambda_K)~dx-|\mathbb{S}^n|e^{\frac{p\left(\widetilde{E}(L)-\widetilde{E}(K)\right)}{|\mathbb{S}^n|}}\right)
		\nonumber
		\\
		\geq&\left(\frac{2n+2+p}{2}-{\eta}\right)\delta_2^2(\widetilde{K},L)
		\\
		=&N_0^{-2}\left(\frac{2n+2+p}{2}-{\eta}\right)\delta_2^2( K,N_0L),
	\end{align*}
	then above inequality yields \eqref{equ:dual p-M-L2-q=0}.
	
	If $p = 0$, similarly substituting \eqref{equ:equ:DU>(2n+1)u} into \eqref{equ:deficit=p=0} yields
	\begin{align*}
		&\int_{\mathbb{S}^n}\log\left(\frac{u_L}{u_K}\right)u_K\rho_K^{-n-1}\sigma_n(\lambda_K)~dx+\widetilde{E}(K)-\widetilde{E}(L)\nonumber
		\\
		=&\int_{\mathbb{S}^n}\left(\log(\frac{u_L}{u_{\widetilde{K}}})\right)u_{\widetilde{K}}\rho_{\widetilde{K}}^{-n-1}\sigma_n(\lambda_{\widetilde{K}})~dx
		\\
		\overset{\eqref{equ:deficit=p=0}}{=}&\frac{1}{2}\Vert\nabla\bar{w}-\nabla \widetilde{w}\Vert_{L^2}^2+\mathscr{E}
		\\
		\overset{\eqref{equ:equ:DU>(2n+1)u}}{\geq}& \left((n+1)-\eta\right)\delta_2^2(\widetilde{K},L)
		\\
		\geq& N_0^{-2}\left((n+1)-\eta\right)\delta_2^2(K,N_0L).
	\end{align*} 
	That is \eqref{equ:dual p-M-L2-p=0-q=0}.
\end{proof}

\subsection{In terms of relative asymmetry index}

In this subsection, we first present Theorem \ref{thm:p-M-L2-F-Vk=Vk} below, which is equivalent to Corollary \ref{thm:dual p-M-al}. We then prove Corollary \ref{thm:dual p-M-al} (in the unnormalized case). 
\begin{theorem}
	\label{thm:p-M-L2-F-Vk=Vk}
	Let $K, L \in \mathcal{K}^{n+1}$ be $C^3$ convex bodies and $p,q$ be real numbers with  $p > q - 2n - 2$. 
	For all $\eta>0$, there exists $\varepsilon>0$  such that  if 
	$K,L\in\mathcal{N}_{\varepsilon}$, $K,L$ satisfy  $S_p(K)=S_p(L)$, and either $\widetilde{W}_q(K)=\widetilde{W}_q(L)$ (for $q \neq 0$) or $\widetilde{E}(K)=\widetilde{E}(L)$ (when $q = 0$), then for $q\neq 0, p=0$, we have
	\begin{align}\label{dual p-M-al-p=0-q!=0-W=W}
		\int_{\mathbb{S}^n}\left(\log(\frac{u_L}{u_K})\right)u_K\rho_K^{q-n-1}\sigma_n(\lambda_K) ~dx
		\geq \widetilde{W}_q(K)\left(c(n,0,q)-\eta\right)\alpha(K,L)^2,
	\end{align}
	and for $q\neq0, p\neq 0$, we have
	\begin{align}\label{dual p-M-al-p!=0-q!=0-W=W}
		\frac{1}{p}\left(\int_{\mathbb{S}^n}u_L^pu_K^{1-p}\rho_K^{q-n-1}\sigma_n(\lambda_K)~dx-\widetilde{W}_q(K)\right)
		\geq \widetilde{W}_q(K)\left(c(n,p,q)-\eta\right)\alpha(K,L)^2,
	\end{align}
	and  for $q=0, p\neq 0$, we have
	\begin{align}\label{equ:dual p-M-al-q=0-W=W}
		&\frac{1}{p}\left(\int_{\mathbb{S}^n}u_L^pu_K^{1-p}\rho_K^{-n-1}\sigma_n(\lambda_K)~dx-|\mathbb{S}^n|\right)
		\nonumber
		\\
		\geq& |\mathbb{S}^n|\left(c(n,p,0)-{\eta}\right)\alpha( K,L)^2,
	\end{align}
	and for $p=0,q=0$, we have 
	\begin{align}
		\label{equ:dual p-M-al-p=0-q=0-W=W}
		\int_{\mathbb{S}^n}\log\left(\frac{u_L}{u_K}\right)u_K\rho_K^{-n-1}\sigma_n(\lambda_K)~dx\nonumber
		\\
		\geq |\mathbb{S}^n|\left(c(n,0,0)-\eta\right)\alpha(K,L)^2.
	\end{align}
	Here $N,N_0$ as defined in Corollary \ref{thm:dual p-M-L2}.
\end{theorem} 

\begin{corollary}
	\label{thm:dual p-M-al} 
	Suppose $K, L \in \mathcal{K}^{n+1}$ are $C^3$ convex bodies, $p$ and $q$ are real numbers with  $p > q - 2n - 2$.
	For all $\eta>0$, there exists $\varepsilon>0$  such that if 
	$K,L\in\mathcal{N}_{\varepsilon}$, and $K,L$  satisfy  $S_p(K)=S_p(L)={o}$, then for $p\neq 0, q\neq 0$, we have
	\begin{align}\label{equ:dual p-M-al}
		&\frac{1}{p}\left(\int_{\mathbb{S}^n}u_L^pu_K^{1-p}\rho_K^{q-n-1}\sigma_n(\lambda_K)~dx-\widetilde{W}_q(K)^{1-\frac{p}{q}}\widetilde{W}_q(L)^{\frac{p}{q}}\right)\nonumber
		\\
		\geq& N^{-p}\widetilde{W}_q(K)\left(c(n,p,q)-{\eta}\right)\alpha( K,L)^2,
	\end{align}
	and for $p=0, q\neq 0$, we have
	\begin{align}\label{equ:dual p-M-al-p=0}
		&\int_{\mathbb{S}^n}\left(\log(\frac{u_L}{u_K})\right)u_K\rho_K^{q-n-1}\sigma_n(\lambda_K) ~dx
		-\frac{1}{q}\widetilde{W}_q(K)\log\left(\frac{\widetilde{W}_q(L)}{\widetilde{W}_q(K)}\right)\nonumber
		\\
		\geq& \widetilde{W}_q(K)\left(c(n,0,q)-\eta\right)\alpha( K,L)^2,
	\end{align}
	and for $p\neq 0, q=0$, 
	we have
	\begin{align}\label{equ:dual p-M-al-q=0}
		&\frac{1}{p}\left(\int_{\mathbb{S}^n}u_L^pu_K^{1-p}\rho_K^{-n-1}\sigma_n(\lambda_K)~dx-|\mathbb{S}^n|e^{\frac{p\left(\widetilde{E}(L)-\widetilde{E}(K)\right)}{|\mathbb{S}^n|}}\right)
		\nonumber
		\\
		\geq& |\mathbb{S}^n|N_0^{-p}\left(c(n,p,0)-{\eta}\right)\alpha( K,L)^2,
	\end{align}
	and for $p=0,q=0$, we have 
	\begin{align}
		\label{equ:dual p-M-al-p=0-q=0}
		\int_{\mathbb{S}^n}\log\left(\frac{u_L}{u_K}\right)u_K\rho_K^{-n-1}\sigma_n(\lambda_K)~dx+\widetilde{E}(K)-\widetilde{E}(L)\nonumber
		\\
		\geq |\mathbb{S}^n|\left(c(n,0,0)-\eta\right)\alpha(K, L)^2.
	\end{align}
\end{corollary} 

The proofs of Theorem \ref{thm:p-M-L2-F-Vk=Vk} and Corollary \ref{thm:dual p-M-al} rely on an estimate for the relative asymmetry index $\alpha(K,L)$, defined in \eqref{equ:def-asymmetry}. Before establishing this estimate, we need  following  lemma.

\begin{lemma}[\cite{Ding-Gao-Li-BM1}*{Lemma 4.5}]
	Suppose $K, L \in \mathcal{K}_o^{n+1}$ are $C^2$ convex bodies, with support function $u_K$ and $u_L$, respectively.  Then we have
	\begin{align*}
		\left|K\Delta L\right|
		\leq \int_0^1\int_{\mathbb{S}^n}
		\left|u_K-u_L
		\right|
		\cdot \det\left(\nabla^2u_t+u_t I\right)
		dx dt
		.
	\end{align*}
	Here $u_t$ is the support function of $(1-t)K+tL$ ($t\in[0,1]$).
	Moreover, if $K, L \in\mathcal{N}_{\varepsilon}$, we have
	\begin{align}
		\left|K\Delta L\right|
		\leq\left(1+O(\varepsilon)\right)\cdot |\mathbb{S}^n|^\frac12\cdot \Vert u_K-u_L\Vert_{L^2}
		.\label{equ:pf-lemma-alpha>-4}
	\end{align}
\end{lemma}
	\begin{lemma}
	\label{lemma:asymmetry}
	Suppose $K, L \in \mathcal{N}_{\varepsilon}$ are $C^3$ convex bodies.  If $\widetilde{W}_q(K)=\widetilde{W}_q(L)$ (for $q \neq 0$) or $\widetilde{E}(K)=\widetilde{E}(L)$, then we have
	\begin{align}
		\label{equ:lemma-asymmetry}
		\alpha(K,L)^2\leq \frac{(n+1)^2}{|\mathbb{S}^n|}\cdot
		\|w-\bar{w}\|_{L^2}^2
		+O(\varepsilon)\cdot
		\|w-\bar{w}\|_{L^2}^2.
	\end{align}
\end{lemma}
\begin{proof}
	
	Since $K,L\in \mathcal{N}_{\varepsilon}$, we have $K_t:=(1-t)K+tL\in\mathcal{N}_{\varepsilon}$ with $t\in[0,1]$ (from Lemma \ref{lemma:C2-1-p-sum-C2-1} (i)). Let $u_t$ and $\rho_t$ denote the support function and radial function of $K_t$, respectively. Then from \eqref{equ:rho^} and \eqref{sigma n=1+o}, we have
	\begin{align}
		\label{equ:pf-lemma-alpha>-5}
		\left(1-\rho_{t}^{q-n-1}\right)\sigma_n\left(\nabla^2u_t+u_t I\right)=O(\varepsilon).
	\end{align}
	If  $\widetilde{W}_q(K)=\widetilde{W}_q(L)$, from \eqref{equ:pf-Lemma 3.2},  for $q\neq 0$, we have:
	\begin{align}
		\label{equ:pf-lemma-alpha>-6}
		0=\widetilde{W}_q(L)-\widetilde{W}_q(K)=q\int_0^1\int_{\mathbb{S}^n}v\rho_{t}^{q-n-1}\sigma_n(\tau[u_t])dxdt.
	\end{align}
	If $\widetilde{E}(K)=\widetilde{E}(L)$, from \eqref{equ:pf-Lemma 3.2-0},
	we have
	\begin{align}
		\label{equ:pf-lemma-alpha>-6-0}
		0=\widetilde{E}(L)-\widetilde{E}(K)=\int_0^1\int_{\mathbb{S}^n}v\rho_{t}^{-n-1}\sigma_n(\tau[u_t])dxdt.
	\end{align}
	Since $\widetilde{W}_{n+1}(K)=(n+1)|K|$, then we have
	\begin{align}
		\label{equ:pf-lemma-alpha>-7}
		\big||L|-|K|\big|=&\left|\int_0^1\int_{\mathbb{S}^n}v\sigma_n(\tau[u_t])dxdt\right|\nonumber
		\\
		\overset{\eqref{equ:pf-lemma-alpha>-6},\eqref{equ:pf-lemma-alpha>-6-0}}{=}&\left|\int_0^1\int_{\mathbb{S}^n}v\left(\sigma_n(\tau[u_t])-\rho_{t}^{q-n-1}\sigma_n(\tau[u_t])\right)dxdt\right|\nonumber
		\\
		\overset{\eqref{equ:pf-lemma-alpha>-5}}{\leq}&O(\varepsilon)\Vert v\Vert_{L^2},
	\end{align}
	where we used the H\"older inequality.
	Let $\W_{K}=rL$, where $r>0$ is a constant chosen such that $|K|=|\W_{K}|=r^{n+1}|L|$, then
	\begin{align}\label{equ:pf-lemma-alpha>-8}
		r^{n+1}=\frac{|K|}{|L|}.
	\end{align}	
	Combining  with  $|L|=\frac{|\mathbb{S}^n|}{(n+1)}+O(\varepsilon)$, and $|K|=\frac{|\mathbb{S}^n|}{(n+1)}+O(\varepsilon)$ we have $r=1+O(\varepsilon)$ and $rL\in\mathcal{N}_{O(\varepsilon)}$. From \eqref{equ:pf-lemma-alpha>-8}, we have
	\begin{align*}
		r-1=\frac{|K|-|L|}{|L|\sum_{j=0}^nr^j},
	\end{align*}
	which implies
	\begin{align}
		\label{equ:pf-lemma-alpha>-9}
		\left|
		r-1\right|\leq O(\varepsilon)\Vert v\Vert_{L^2}.
	\end{align}
	Then we have
	\begin{align*}
		\Vert ru_L-u_K\Vert_{L^2}=&\Vert (r-1)u_L+v\Vert_{L^2}
		\\
		\leq &|r-1|\cdot \Vert u_L\Vert_{L^2}+\Vert v\Vert_{L^2}
		\\
		\leq& \left(1+O(\varepsilon)\right)\Vert v\Vert_{L^2}.
	\end{align*}
	Putting it into \eqref{equ:pf-lemma-alpha>-4} yields
	\begin{align*}
		|K\Delta\W_{K}|=|K\Delta(rL)|\leq \left(1+O(\varepsilon)\right)\cdot |\mathbb{S}^n|^\frac12\cdot \Vert v\Vert_{L^2}.
	\end{align*}
	Since $|\W_K|={|K|}=\frac{|\mathbb{S}^n|}{n+1}(1+O(\varepsilon))$, 
	we have
	\begin{align*}
		\frac{|K\Delta\W_{K}|^2}{|\W_{K}|^2}
		\leq \left(\frac{(n+1)^2}{|\mathbb{S}^n|}
		+O(\varepsilon)\right)
		\|w-\bar{w}\|_{L^2}^2.
	\end{align*}
	On the other hand,
	\begin{align*}
		\alpha(K,L)=\min _{y \in \mathbb{R}^{n+1}}\left\{\frac{|K \Delta(r L+y)|}{|K|}:|r L|=|K|\right\}\leq \frac{|K \Delta \W_{K}|}{|\W_{K}|}.
	\end{align*}
	Thus we obtain
	\begin{align*}
		\alpha(K,L)^2\leq \left(\frac{(n+1)^2}{|\mathbb{S}^n|}
		+O(\varepsilon)\right)
		\|w-\bar{w}\|_{L^2}^2.
	\end{align*}
\end{proof}




Now we can prove the Theorem \ref{thm:p-M-L2-F-Vk=Vk} as follows:
\begin{proof}[\textbf{Proof of Theorem \ref{thm:p-M-L2-F-Vk=Vk}}]
	From \eqref{equ:Wq=} and \eqref{equ:E(K)=}, we know that for $K \in \mathcal{N}_\varepsilon$, $$\widetilde{W}_q(K)=|\mathbb{S}^n|\left(1+O(\varepsilon)\right),\,\text{for}\ q\neq 0
	.$$	
	We Denote 
		\begin{align}\label{equ:W=W E=E main neq}
		\Delta:=	\frac{p-q}{2}\Vert\bar{w}-w\Vert_{L^2}^2
		+\frac{1}{2}\Vert\nabla\bar{w}-\nabla w\Vert_{L^2}^2+\mathscr{E}.	
	\end{align}
	By Lemma \ref{lemma:I_k-w-barw-O(e)}, it suffices to prove the following inequality for any $p>q -2n-2$, 
	\begin{align*}		
	\Delta	\geq |\mathbb{S}^n|c(n,p,q)\left(1+O(\varepsilon)\right)\alpha(K,L)^2.
	\end{align*}
	We now prove this inequality as follows.
	
	 Putting \eqref{equ:equ:DU>(2n+1)u} and \eqref{equ:lemma-asymmetry} into \eqref{equ:W=W E=E main neq}, we obtain
	\begin{align}
		\label{equ:Delta>w-w>}
		\Delta\overset{\eqref{equ:equ:DU>(2n+1)u}}{\geq}& \left(\frac{p-q}{2}+n+1\right)
		\left(1+O(\varepsilon)\right)\Vert\bar{w}- w\Vert_{L^2}^2
		\nonumber
		\\
		\overset{\eqref{equ:lemma-asymmetry}}{\geq}&
		|\mathbb{S}^n|\frac{(2n+2+p-q)}{2(n+1)^2}
		\left(1+O(\varepsilon)\right)\alpha(K,L)^2.
	\end{align}

	Then we complete the proof. 
\end{proof}

\begin{proof}[\textbf{Proof of Corollary \ref{thm:dual p-M-al}}] 	
	We only verify the case $p = 0$, $q \neq 0$; the other cases follow similarly from Theorem \ref{thm:p-M-L2-F-Vk=Vk}.	For $K$ and $L$, there exists an constant $N>0$, such that $\widetilde{K}:=\frac{1}{N}K$ satisfies $\widetilde{W}_q(\widetilde{K})=\widetilde{W}_q(L)$ with $q\neq 0$. Then we have
	$$\left(\frac{1}{N}\right)^{q}\widetilde{W}_q(K)=\widetilde{W}_q(L),$$
	which implies 
	\begin{align*}
		N=\left(\frac{\widetilde{W}_q(K)}{\widetilde{W}_q(L)}\right)^{\frac{1}{q}}.
	\end{align*}
	From \eqref{equ:Wq=}, we have $\widetilde{W}_q(K)= |\mathbb{S}^{n}|+O(\varepsilon)$, and $\widetilde{W}_q(L)=|\mathbb{S}^{n}|+O(\varepsilon)$, then $N=1+O(\varepsilon)$. Denote $u_{\widetilde{K}}=\frac{u_K}{N}:=1+\tilde{w}$ is the support function of $\widetilde{K}$, then we have $\Vert \tilde{w}\Vert_{C^2}=O(\varepsilon)$.  We can check that  pair $(\widetilde{K},L)$ satisfies $S_p(\widetilde{K})=S_p(L)={o}$.
	
	By $K=N\widetilde{K}$ and Theorem \ref{thm:p-M-L2-F-Vk=Vk}, we have  
	\begin{align*}
		&\int_{\mathbb{S}^n}\left(\log(\frac{u_L}{u_K})\right)u_K\rho_K^{q-n-1}\sigma_n(\lambda_K) ~dx
		-\frac{1}{q}\widetilde{W}_q(K)\log\frac{\widetilde{W}_q(L)}{\widetilde{W}_q(K)}
		\\
		=&\int_{\mathbb{S}^n}\left(\log(\frac{u_L}{Nu_{\widetilde{K}}})\right)Nu_{\widetilde{K}}N^{q-1}\rho_{\widetilde{K}}^{q-n-1}\sigma_n(\lambda_{\widetilde{K}}) ~dx-\frac{1}{q}N^{q}\widetilde{W}_q(\widetilde{K})\log\frac{\widetilde{W}_q(L)}{N^{q}\widetilde{W}_q(\widetilde{K})}
		\\
		\geq& N^{q}\left(c(n,0,q)-\eta\right)\alpha(L,\widetilde{K})^2\widetilde{W}_q(\widetilde{K})
		\\
		=& {\widetilde{W}_q(K)}\left(c(n,0,q)-\eta\right)\alpha(L,\widetilde{K})^2
		.
	\end{align*}
	Since 
	\begin{align*}
		\alpha(L,\widetilde{K})=&\min _{y \in \mathbb{R}^{n+1}}\left\{\frac{|\widetilde{K} \Delta(r L+y)|}{|\widetilde{K}|}:|r L|=|\widetilde{K}|\right\}
		\\
		=&\min _{y \in \mathbb{R}^{n+1}}\left\{\frac{|(\frac{1}{N}K) \Delta(\frac{1}{N}rN L+y)|}{|{\frac{1}{N}K}|}:|r L|=|\frac{1}{N}K|\right\}
		\\
		=&\min _{y \in \mathbb{R}^{n+1}}\left\{\frac{|K \Delta(rN L+y)|}{|{K}|}:|rN L|=|K|\right\}
		\\
		=&\alpha(K,L).
	\end{align*} 
	Then we obtain the inequality \eqref{equ:dual p-M-al-p=0}.
	
\end{proof}

\section{Stability of dual $L^p$-Brunn--Minkowski inequality}\label{sec 5}
In this section, we prove the stability of the dual  $L^p$-Brunn--Minkowski inequality stated in Theorems \ref{thm:dual p-BM-L2-Vk!=Vk} and \ref{thm:dual p-BM-alpha-Vk!=Vk}. We first derive the normalized versions (i.e.,  requiring $\widetilde{W}_q(K)=\widetilde{W}_q(L)$ or $\widetilde{E}(K)=\widetilde{E}(L)$) as given in Theorems \ref{thm:dual p-BM-L2-Vk=Vk}. Before presenting the proofs, We first introduce some necessary lemmas, which are similar to Lemma 5.1  in \cite[Section 5]{Ding-Gao-Li-BM1}.

\begin{lemma}\label{lemma:C2-1-p-sum-C2-1}
	Assume $p,q\in\mathbb{R}$.	Given a $\eta>0$, there exists an $\varepsilon$ such that, for every $K,L\in\mathcal{N}_{\varepsilon}$ and $t\in[0,1]$, the following hold
	\begin{itemize}
		\item [(i)] $(1-t) \cdot K+_p t \cdot L \in \mathcal{N}_{\eta}$,
		\item [(ii)] $u_{(1-t) \cdot K+_p t \cdot L}=\left((1-t) u_K^p+t u_L^p\right)^{1 / p}$ for $p\neq 0$,
		\item [(iii)] $u_{(1-t) \cdot K+_0 t \cdot L}= u_K^{1-t} u_L^{t}$, 
		\item [(iv)] $\left(\frac{\widetilde{W}_q(\mathbb{B}^{n+1})}{\widetilde{W}_q(K)}\right)^{\frac{1}{q}} K \in \mathcal{N}_{\eta}$ for $q\neq 0$, and $e^{-\frac{\widetilde{E}(K)}{|\mathbb{S}^n|}}K\in \mathcal{N}_{\eta}$.
	\end{itemize}
\end{lemma}
\begin{proof}
	Parts (i)$\thicksim$(iii) have been proved in \cite[Lemma 5.1]{Ding-Gao-Li-BM1}, while part (iv) follows from the facts that $\widetilde{W}_q(K)=|\mathbb{S}^{n}|+O(\varepsilon)$ for $q\neq 0$ and $\widetilde{E}(K)=O(\varepsilon)$.
\end{proof}

Using Corollary \ref{thm:dual p-M-L2} and \ref{thm:dual p-M-al}, we first obtain the following dual $L^p$-Minkowski inequalities for the case either $\widetilde{W}_q(K)=\widetilde{W}_q(L)$  where $q \neq 0$ or $\widetilde{E}(K)=\widetilde{E}(L)$.
\begin{theorem}
	\label{thm:dual p-BM-L2-Vk=Vk}
	Suppose $K, L \in \mathcal{K}^{n+1}$ are $C^3$ convex bodies, $p$ and $q$ are real numbers with  $p > q - 2n - 2$, and $t \in [0,1]$. Assume that 
	either $\widetilde{W}_q(K)=\widetilde{W}_q(L)$ for $q\neq0$ or $\widetilde{E}(K)=\widetilde{E}(L)$ when $q=0$. For all $\eta>0$, there exists $\varepsilon>0$  such that if $K,L\in\mathcal{N}_{\varepsilon}$ and $K,L$ satisfy  $S_p(K)=S_p(L)={o}$, then 
	for $p=0, q\neq 0$, we have 
	\begin{align}
		\label{equ:stab L2-BM-Wq=Wq-p=0}
		&\frac{1}{q}\left(\log\widetilde{W}_q\((1-t)\cdot K +_p t\cdot L\)-\log\widetilde{W}_q(K)\right)\nonumber
		\\
		\geq &t(1-t)\widetilde{W}_q(K)^{-1}\left(\frac{2n+2-q}{2}-\eta\right)\delta_2^2(K,L),
	\end{align}
	\begin{align}
		\label{equ:stab a-BM-Wq=Wq-p=0}
		&\frac{1}{q}\left(\log\widetilde{W}_q\((1-t)\cdot K +_p t\cdot L\)-\log\widetilde{W}_q(K)\right)\nonumber
		\\
		\geq &t(1-t)\left(c(n,0,q)-\eta\right)\alpha(K,L)^2,
	\end{align}	
	and for $p\neq 0$ and $q\neq0$, we have
	\begin{align}
		\label{equ:stab L2-BM-Wq=Wq}
		&\frac{1}{p}\left(\widetilde{W}_q\((1-t)\cdot K +_p t\cdot L\)^\frac{p}{q}-\widetilde{W}_q(K)^{\frac{p}{q}}\right)\nonumber
		\\
		\geq &t(1-t)\widetilde{W}_q(K)^{\frac{p-q}{q}}\left(\frac{2n+2+p-q}{2}-\eta\right)\delta_2^2(K,L),
	\end{align}
	\begin{align}
		\label{equ:stab a-BM-Wq=Wq}
		&\frac{1}{p}\left(\widetilde{W}_q\((1-t)\cdot K +_p t\cdot L\)^\frac{p}{q}-\widetilde{W}_q(K)^{\frac{p}{q}}\right)\nonumber
		\\
		\geq &t(1-t)\widetilde{W}_q(K)^{\frac{p}{q}}\left(c(n,p,q)-\eta\right)\alpha(K,L)^2,
	\end{align}
	and for $p\neq 0, q=0$, we have
	\begin{align}\label{equ:stab L2-BM-Wq=Wqp!0-q=0}
		&\frac{1}{p}\left(e^{\frac{p\left(\widetilde{E}((1-t)\cdot K +_p t\cdot L)-\widetilde{E}(K)\right)}{|\mathbb{S}^n|}}-1\right)
		\nonumber
		\\
		\geq& t(1-t)|\mathbb{S}^n|^{-1}\left(\frac{2n+2+p}{2}-{\eta}\right)\delta_2^2( K,L),
	\end{align}	
	\begin{align}\label{equ:stab a-BM-Wq=Wqp!0-q=0}
		&\frac{1}{p}\left(e^{\frac{p\left(\widetilde{E}((1-t)\cdot K +_p t\cdot L)-\widetilde{E}(K)\right)}{|\mathbb{S}^n|}}-1\right)
		\nonumber
		\\
		\geq& t(1-t)\left(c(n,p,0)-{\eta}\right)\alpha( K,L)^2,
	\end{align}	
	and for $p=q=0$, we have
	\begin{align}\label{equ:stab L2-BM-Wq=Wq-q=0}
		&\widetilde{E}\((1-t)\cdot K +_p t\cdot L\)-\widetilde{E}(K)\nonumber
		\\
		\geq&t(1-t)
		\left(\frac{2n+2+p}{2}-\eta\right)\delta_2^2(K,L),
	\end{align}
	\begin{align}\label{equ:stab a-BM-Wq=Wq-q=0}
		&\widetilde{E}\((1-t)\cdot K +_p t\cdot L\)-\widetilde{E}(K)\nonumber
		\\
		\geq&t(1-t)
		|\mathbb{S}^n|\left(c(n,p,0)-\eta\right)\alpha(K,L)^2.
	\end{align}
\end{theorem} 
	\begin{proof}[\textbf{Proof of Theorem \ref{thm:dual p-BM-L2-Vk=Vk}}]
		(1).  We consider the case $q\neq 0, p=0$. From Lemma \ref{lemma:C2-1-p-sum-C2-1},  for  sufficiently small $\varepsilon$,  $K,L\in\mathcal{N}_{\varepsilon}$, and $t\in[0,1]$, we have $Q_t:=(1-t)\cdot K+_0t\cdot L\in\mathcal{N}_{O(\varepsilon)}$ is convex body with  support function $u_{Q_t}=u_K^{1-t} u_L^{t}$. Since $K$ and $L$ satisfy $S_p(K)=S_p(L)={o}$, we know  $S_p(Q_t)={o}$. 
		Then  from \eqref{equ:dual p-M-L2-p=0} we have
		\begin{align}
			\label{equ:pf-thm:p-BM-L2-Vk=Vk-4}
			&-t\int_{\mathbb{S}^n}\left(\log(\frac{u_L}{u_K})\right)u_{Q_t}\rho_{Q_t}^{q-n-1}\sigma_n(\lambda_{Q_t}) ~dx\nonumber
			\\
			=&\int_{\mathbb{S}^n}\left(\log(\frac{u_{K}}{u_{Q_t}})\right)u_{{Q_t}}\rho_{Q_t}^{q-n-1}\sigma_n(\lambda_{Q_t}) ~dx
			\nonumber
			\\
			\overset{\eqref{equ:dual p-M-L2-p=0}}{\geq}&\frac{1}{q}\widetilde{W}_q(Q_t)\log\frac{\widetilde{W}_q(K)}{\widetilde{W}_q(Q_t)} 
			+ N_1^{q-2}\left(\frac{2n+2-q}{2}-\eta\right)\delta_2^2(N_1 K,Q_t).
		\end{align}
		And
		\begin{align}
			\label{equ:pf-thm:p-BM-L2-Vk=Vk-5}
			&(1-t)\int_{\mathbb{S}^n}\left(\log(\frac{u_L}{u_K})\right)u_{Q_t}\rho_{Q_t}^{q-n-1}\sigma_n(\lambda_{Q_t}) ~dx
			\nonumber\\
			=&\int_{\mathbb{S}^n}\left(\log(\frac{u_L}{u_{Q_t}})\right)u_{{Q_t}}\rho_{Q_t}^{q-n-1}\sigma_n(\lambda_{Q_t}) ~dx
			\nonumber
			\\
			\geq&\frac{1}{q}\widetilde{W}_q(Q_t)\log\frac{\widetilde{W}_q(L)}{\widetilde{W}_q(Q_t)} 
			+ N_2^{q-2}\left(\frac{2n+2-q}{2}-\eta\right)\delta_2^2(N_2 L,Q_t).
		\end{align}
		Here $N_1=\left(\frac{\widetilde{W}_q(Q_t)}{\widetilde{W}_q(K)}\right)^{\frac{1}{q}}=\left(\frac{\widetilde{W}_q(Q_t)}{\widetilde{W}_q(L)}\right)^{\frac{1}{q}}=N_2$, since $\widetilde{W}_q(K)=\widetilde{W}_q(L)$. 	
		Multiplying equation \eqref{equ:pf-thm:p-BM-L2-Vk=Vk-4} by $(1-t)$
		and adding it to equation \eqref{equ:pf-thm:p-BM-L2-Vk=Vk-5} multiplied by $t$ yields:
		\begin{align}
			\label{equ:pf-thm:p-BM-L2-Vk=Vk-3}
			&- \frac{1}{q}\widetilde{W}_q(Q_t)\log\frac{\widetilde{W}_q(K)}{\widetilde{W}_q(Q_t)} \nonumber
			\\
			\geq & N_1^{q-2}\left(\frac{2n+2-q}{2}-\eta\right)\left((1-t)\delta_2^2(N_1 K,Q_t)+t\delta_2^2(N_1 L,Q_t)\right).
		\end{align}

		On the other hands, Applying the triangle inequality and H\"older's inequality, we obtain 
		\begin{align*}
			\delta_2^2(K,L)\leq& \left(\delta_2(K,Q_t)+\delta_2(Q_t,L)\right)^2
			\nonumber
			\\
			=&\delta_2^2(K,Q_t)+\delta_2^2(Q_t,L)+2\delta_2(K,Q_t)\cdot \delta_2(Q_t,L)
			\nonumber
			\\
			\leq &(1+\frac{1-t}{t})\delta_2^2(K,Q_t)+(1+\frac{t}{1-t})\delta_2^2(Q_t,L)
			\nonumber
			\\
			=&
			\frac{1}{t(1-t)}\left((1-t)\delta_2^2(K,Q_t)+t\delta_2^2(Q_t,L)\right),
		\end{align*}
		where $t\in(0,1)$.  This means for $t\in[0,1]$ and any $N_1>0$, we have
		\begin{align}
			&\label{equ:pf-thm:p-BM-L2-Vk=Vk-2}
			(1-t)\delta_2^2(N_1K,Q_t)+t\delta_2^2(Q_t,N_1L)\nonumber
			\\
			\geq& t(1-t)\delta_2^2(N_1K,N_1L)\nonumber
			\\
			=&N_1^2t(1-t)\delta_2^2(K,L).
		\end{align}
		Putting \eqref{equ:pf-thm:p-BM-L2-Vk=Vk-2} into \eqref{equ:pf-thm:p-BM-L2-Vk=Vk-3} yields \eqref{equ:stab L2-BM-Wq=Wq-p=0}. 
		\\
		(2). We consider the case  $q\neq 0,p\neq 0$. From Lemma \ref{lemma:C2-1-p-sum-C2-1},  for  sufficiently small $\varepsilon$,  $K,L\in\mathcal{N}_{\varepsilon}$, and $t\in[0,1]$, we have $Q_t:=(1-t)\cdot K+_pt\cdot L\in\mathcal{N}_{O(\varepsilon)}$ is convex body with  support function $u_{Q_t}=\left((1-t)u_K^p+tu_L^p\right)^{\frac{1}{p}}$.  Since $S_p(K)=S_p(L)={o}$, we have   $S_p(Q_t)={o}$.
		Then 
		\begin{align}\label{equ:pf-thm:p-BM-L2-Vk=Vq-2}
			&\frac{1}{p}\widetilde{W}_q(Q_t)=\frac{1}{p}\int_{\mathbb{S}^n}u_{Q_t}\rho_{Q_t}^{q-n-1}\sigma_n(\lambda_{Q_t})~dx\nonumber
			\\
			=&\frac{1-t}{p}\int_{\mathbb{S}^n}u_{K}^pu_{Q_t}^{1-p}\rho_{Q_t}^{q-n-1}\sigma_n(\lambda_{Q_t})~dx
			+\frac{t}{p}\int_{\mathbb{S}^n}u_{L}^pu_{Q_t}^{1-p}\rho_{Q_t}^{q-n-1}\sigma_n(\lambda_{Q_t})~dx\nonumber
			\\
			\overset{\eqref{equ:dual p-M-L2}}{\geq}
			&\frac{1-t}{p}\widetilde{W}_q(Q_t)^{\frac{q-p}{q}}\widetilde{W}_q(K)^{\frac{p}{q}}+\frac{t}{p}\widetilde{W}_q(Q_t)^{\frac{q-p}{q}}\widetilde{W}_q(L)^{\frac{p}{q}}
			\nonumber
			\\
			&+{(1-t)N_1^{q-2-p}}\left(\frac{2n+2+p-q}{2}+O(\varepsilon)\right)\delta_2^2(N_1K,Q_t)\nonumber
			\\
			&+{tN_2^{q-2-p}}\left(\frac{2n+2+p-q}{2}+O(\varepsilon)\right)\delta_2^2(Q_t,N_2L)\nonumber
			\\
			=&\frac{1}{p}\widetilde{W}_q(Q_t)^{\frac{q-p}{q}}\left((1-t)\widetilde{W}_q(K)^{\frac{p}{q}}+t\widetilde{W}_q(L)^{\frac{p}{q}}\right)
			\nonumber
			\\
			&+{N_1^{q-2-p}}\left(\frac{2n+2+p-q}{2}+O(\varepsilon)\right)\left((1-t)\delta_2^2(N_1K,Q_t)+t\delta_2^2(Q_t,N_1L)\right).
		\end{align}
		Here $N_1=\left(\frac{\widetilde{W}_q(Q_t)}{\widetilde{W}_q(K)}\right)^{\frac{1}{q}}=\left(\frac{\widetilde{W}_q(Q_t)}{\widetilde{W}_q(L)}\right)^{\frac{1}{q}}=N_2$, since $\widetilde{W}_q(K)=\widetilde{W}_q(L)$.
		
		Putting \eqref{equ:pf-thm:p-BM-L2-Vk=Vk-2} into \eqref{equ:pf-thm:p-BM-L2-Vk=Vq-2} yields \eqref{equ:stab L2-BM-Wq=Wq}.
		\\
		(3). 	For case $q=0,p\neq 0$, we denote $Q_t=(1-t)\cdot K+_p t\cdot L$. From \eqref{equ:dual p-M-L2-q=0} in Corollary \ref{thm:dual p-M-L2}, we have
		\begin{align}\label{equ:pf-thm:p-BM-L2-Vk=Vk-7}
			&\frac{1}{p}\left(\int_{\mathbb{S}^n}u_K^pu_{Q_t}^{1-p}\rho_{Q_t}^{-n-1}\sigma_n(\lambda_{Q_t})~dx-|\mathbb{S}^n|\widetilde{N_0}^{-p}\right)
			\nonumber
			\\
			\geq& \widetilde{N_0}^{-p-2}\left(\frac{2n+2+p}{2}-{\eta}\right)\delta_2^2( {Q_t},\widetilde{N_0}K),
		\end{align}
		and 
		\begin{align}\label{equ:pf-thm:p-BM-L2-Vk=Vk-8}
			&\frac{1}{p}\left(\int_{\mathbb{S}^n}u_L^pu_{Q_t}^{1-p}\rho_{Q_t}^{-n-1}\sigma_n(\lambda_{Q_t})~dx-|\mathbb{S}^n|\widetilde{N_0}^{-p}\right)
			\nonumber
			\\
			\geq& \widetilde{N_0}^{-p-2}\left(\frac{2n+2+p}{2}-{\eta}\right)\delta_2^2( {Q_t},\widetilde{N_0}L),
		\end{align}
		with $\widetilde{N_0}=e^{\frac{\widetilde{E}(Q_t)-\widetilde{E}(L)}{|\mathbb{S}^n|}}$. 
		Multiplying equation \eqref{equ:pf-thm:p-BM-L2-Vk=Vk-7} by $(1-t)$
		and adding it to equation \eqref{equ:pf-thm:p-BM-L2-Vk=Vk-8} multiplied by $t$ yields:
		\begin{align}\label{equ:pf-thm:p-BM-L2-Vk=Vk-9}
			&\frac{1}{p}\left(\int_{\mathbb{S}^n}u_{Q_t}\rho_{Q_t}^{-n-1}\sigma_n(\lambda_{Q_t})~dx-|\mathbb{S}^n|\widetilde{N_0}^{-p}\right)
			\nonumber
			\\
			\geq& \widetilde{N_0}^{-p-2}\left(\frac{2n+2+p}{2}-{\eta}\right)\left((1-t)\delta_2^2( {Q_t},\widetilde{N_0}K)+t\delta_2^2( {Q_t},\widetilde{N_0}L)\right).
		\end{align}
	Substituting the identity $\int_{\mathbb{S}^n}u_{Q_t}\rho_{Q_t}^{-n-1}\sigma_n(\lambda_{Q_t})~dx=|\mathbb{S}^n|$ together with \eqref{equ:pf-thm:p-BM-L2-Vk=Vk-2} into \eqref{equ:pf-thm:p-BM-L2-Vk=Vk-9} yields
	\begin{align}\label{equ:pf-thm:p-BM-L2-Vk=Vk-15}
		&\frac{1}{p}\left(1-e^{-p\frac{\widetilde{E}(Q_t)-\widetilde{E}(L)}{|\mathbb{S}^n|}}\right)
		\nonumber
		\\
		\geq&\widetilde{N_0}^{-p} t(1-t)|\mathbb{S}^n|^{-1}\left(\frac{2n+2+p}{2}-{\eta}\right)\delta_2^2( K,L).
	\end{align}	
	 This is inequality \eqref{equ:stab L2-BM-Wq=Wqp!0-q=0}.

		(4). For $q=p=0$, we also denote $Q_t=(1-t)\cdot K+_0 t\cdot L$. From \eqref{equ:dual p-M-L2-p=0-q=0} in Corollary \ref{thm:dual p-M-L2}, we have
		\begin{align}
			\label{equ:pf-thm:p-BM-L2-Vk=Vk-10}
			&\int_{\mathbb{S}^n}\log\left(\frac{u_K}{u_{Q_t}}\right)u_{Q_t}\rho_{Q_t}^{-n-1}\sigma_n(\lambda_{Q_t})~dx+\widetilde{E}({Q_t})-\widetilde{E}(K)\nonumber
			\\
			\geq &\widetilde{N_0}^{-2}\left((n+1)-\eta\right)\delta_2^2({Q_t},\widetilde{N_0}K).
		\end{align}
		and \begin{align}
			\label{equ:pf-thm:p-BM-L2-Vk=Vk-11}
			&\int_{\mathbb{S}^n}\log\left(\frac{u_L}{u_{Q_t}}\right)u_{Q_t}\rho_{Q_t}^{-n-1}\sigma_n(\lambda_{Q_t})~dx+\widetilde{E}({Q_t})-\widetilde{E}(L)\nonumber
			\\
			\geq &\widetilde{N_0}^{-2}\left((n+1)-\eta\right)\delta_2^2({Q_t},\widetilde{N_0}L).
		\end{align}
		Multiplying equation \eqref{equ:pf-thm:p-BM-L2-Vk=Vk-10} by $(1-t)$
		and adding it to equation \eqref{equ:pf-thm:p-BM-L2-Vk=Vk-11} multiplied by $t$ yields \eqref{equ:stab L2-BM-Wq=Wq-q=0}. 

		The proof of \eqref{equ:stab a-BM-Wq=Wq-p=0}, \eqref{equ:stab a-BM-Wq=Wq}, \eqref{equ:stab a-BM-Wq=Wqp!0-q=0},  and \eqref{equ:stab a-BM-Wq=Wq-q=0} are analogous to that of \eqref{equ:stab L2-BM-Wq=Wq-p=0}, \eqref{equ:stab L2-BM-Wq=Wq}, \eqref{equ:stab L2-BM-Wq=Wqp!0-q=0},  and \eqref{equ:stab L2-BM-Wq=Wq-q=0} above. Moreover, using \cite[eq. (5.14)]{Ding-Gao-Li-BM1}, we have for $t \in [0,1]$:
		\begin{align}
			\label{equ:pf-thm:p-BM-L2-Vk=Vk-2-index}
			(1-t)\alpha(K,Q_t)^2 + t\alpha(L,Q_t)^2 \ge t(1-t) \alpha(K,L)^2.
		\end{align}
		For brevity, we omit the proof of \eqref{equ:stab a-BM-Wq=Wq-p=0}, \eqref{equ:stab a-BM-Wq=Wq}, \eqref{equ:stab a-BM-Wq=Wqp!0-q=0},   and \eqref{equ:stab a-BM-Wq=Wq-q=0}.
	\end{proof}
	
	Then we can derive Theorem \ref{thm:dual p-BM-L2-Vk!=Vk} and Theorem \ref{thm:dual p-BM-alpha-Vk!=Vk} from Theorem \ref{thm:dual p-BM-L2-Vk=Vk}. Since the proof of Theorem \ref{thm:dual p-BM-alpha-Vk!=Vk} is similar to that of Theorem \ref{thm:dual p-BM-L2-Vk!=Vk}, we will only prove the latter.

	\begin{proof}[\textbf{Proof of Theorem \ref{thm:dual p-BM-L2-Vk!=Vk}}]
		For case $q\neq0, p=0$,	we take $\tilde{K}=N_K  K,\tilde{L}=N_L  L$, where $N_K=\left(\frac{\widetilde{W}_q(\mathbb{B}^{n+1})}{\widetilde{W}_q(K)}\right)^{\frac{1}{q}}$ and $N_L=\left(\frac{\widetilde{W}_q(\mathbb{B}^{n+1})}{\widetilde{W}_q(L)}\right)^{\frac{1}{q}}$. Then we know $\widetilde{W}_q(\tilde{K})=\widetilde{W}_q(\tilde{L})=\widetilde{W}_q(\mathbb{B}^{n+1})=|\mathbb{S}^n|$ and $\tilde{K},\tilde{L}$ are also  $C^3$  convex bodies in $\mathcal{N}_{O(\varepsilon)}$ by Lemma \ref{lemma:C2-1-p-sum-C2-1}. 
		Sinve $S_p(K)=S_p(L)={o}$, we have that $\tilde{K}$ and $\tilde{L}$ also satisfy $S_p(\tilde{K})=S_p(\tilde{L})={o}$. 
		Take
		\begin{align*}
			\beta=N_K^{1-t}N_L^t,
		\end{align*} 
		then $(1-t)\cdot \tilde{K}+_0t\cdot \tilde{L}=\beta \left((1-t)\cdot K+_0t\cdot L\right)$  by Lemma \ref{lemma:C2-1-p-sum-C2-1} (iii). 
		Form \eqref{equ:stab L2-BM-Wq=Wq-p=0}, we have
		\begin{align*}
			&\frac{1}{q}\left(\log	\left(\widetilde{W}_q\((1-t)\cdot \tilde{K}+_0t\cdot \tilde{L}\)\right)-\log\left(\widetilde{W}_q(\tilde{K})\right)\right)\nonumber
			\\
			\geq&
			t(1-t)\widetilde{W}_q(\tilde{K})^{-1}\left(\frac{2n+2-q}{2}-\eta\right)\delta_2^2(\tilde{K},\tilde{L}).
		\end{align*}
		This means 
		\begin{align*}
			&\frac{1}{q}\log	\left(\widetilde{W}_q\((1-t)\cdot K+_0t\cdot L\)\right)
			+\frac{1}{q}(1-t)\left(\log\(\widetilde{W}_q(\mathbb{B}^{n+1})\)-\log\( \widetilde{W}_q(K)\)\right)
			\\
			&+\frac{1}{q}t\left(\log\(\widetilde{W}_q(\mathbb{B}^{n+1})\)-\log\( \widetilde{W}_q(L)\)\right)
			-\frac{1}{q}\log\(\widetilde{W}_q(\mathbb{B}^{n+1})\)\nonumber
			\\
			\geq&
			t(1-t)\widetilde{W}_q(\mathbb{B}^{n+1})^{-1}\left(\frac{2n+2-q}{2}-\eta\right)\delta_2^2(\tilde{K},\tilde{L}).
		\end{align*}
		which implies \eqref{equ:stab L2-BM-Wq!=Wq-p=0}.

		For $p\neq 0, q\neq 0$, we take $\tilde{K}=N_K  K,\tilde{L}=N_L  L$, where $N_K=\left(\frac{\widetilde{W}_q(\mathbb{B}^{n+1})}{\widetilde{W}_q(K)}\right)^{\frac{1}{q}}$,  $N_L=\left(\frac{\widetilde{W}_q(\mathbb{B}^{n+1})}{\widetilde{W}_q(L)}\right)^{\frac{1}{q}}$, and take
		\begin{align*}
			\tilde{t}=\frac{t\widetilde{W}_q(L)^{\frac{p}{q}}}{(1-t)\widetilde{W}_q(K)^{\frac{p}{q}}+t\widetilde{W}_q(L)^{\frac{p}{q}}}\in[0,1].
		\end{align*} 
		By Lemma \ref{lemma:C2-1-p-sum-C2-1}(ii), one can  check that 
		\begin{align*}
			(1-\tilde{t})\cdot \tilde{K}+_p \tilde{t}\cdot\tilde{L}=\frac{|\mathbb{S}^n|^{\frac{1}{q}}}{\left({(1-t)\widetilde{W}_q(K)^{\frac{p}{q}}+t\widetilde{W}_q(L)^{\frac{p}{q}}}\right)^{\frac{1}{p}}}\((1-t)\cdot K+_pt\cdot L\).
		\end{align*}
		Form \eqref{equ:stab L2-BM-Wq=Wq}, we have
		\begin{align*}
			&\frac{1}{p}\left(\widetilde{W}_q\((1-\tilde{t})\cdot \tilde{K}+_p\tilde{t}\cdot \tilde{L}\)^{\frac{p}{q}}-\widetilde{W}_q(\tilde{L})^{\frac{p}{q}}\right)\nonumber
			\\
			\geq&\tilde{t}(1-\tilde{t})\widetilde{W}_q(\tilde{K})^{-\frac{q-p}{q}}\left(\frac{2n+2+p-q}{2}-\eta\right)\delta_2^2(\tilde{K},\tilde{L}).
		\end{align*}
		This means 
		\begin{align*}
			&\frac{1}{p}\left(\frac{\widetilde{W}_q\((1-t)\cdot K+_pt\cdot L\)^{\frac{p}{q}}\widetilde{W}_q(\mathbb{B}^{n+1})^{\frac{p}{q}}}{(1-t)\widetilde{W}_q(K)^{\frac{p}{q}}+t\widetilde{W}_q(L)^{\frac{p}{q}}}-\widetilde{W}_q(\mathbb{B}^{n+1})^{\frac{p}{q}}\right)\nonumber
			\\
			\geq& \frac{t(1-t)\widetilde{W}_q(K)^{\frac{p}{q}}\widetilde{W}_q(L)^{\frac{p}{q}}}{\left((1-t)\widetilde{W}_q(K)^{\frac{p}{q}}+t\widetilde{W}_q(L)^{\frac{p}{q}}\right)^2}\widetilde{W}_q(\mathbb{B}^{n+1})^{\frac{p}{q}-1}\left(\frac{2n+2+p-q}{2}-\eta\right)\delta_2^2(\tilde{K},\tilde{L}),
		\end{align*}
		which implies \eqref{equ:stab L2-BM-Wq!=Wq}.

		For $q=0,p\neq 0$, we still  take $\tilde{K}=\widetilde{N_K}  K,\tilde{L}=\widetilde{N_L}  L$ satisfy $\widetilde{E}(\widetilde{K})=\widetilde{E}(\widetilde{L})=\widetilde{E}(\mathbb{B}^{n+1})=0$, then  $$\widetilde{N_K}=e^{-\frac{\widetilde{E}(K)}{|\mathbb{S}^n|}},\quad \widetilde{N_L}=e^{-\frac{\widetilde{E}(L)}{|\mathbb{S}^n|}}.$$
		By Lemma \ref{lemma:C2-1-p-sum-C2-1}, we know $\tilde{K},\tilde{L}\in \mathcal{N}_{O(\varepsilon)}$, and $S_p(\tilde{K})=S_p(\tilde{L})={o}$. 
		We denote
		\begin{align*}
			s=\frac{t\widetilde{N_L}^{-p}}{(1-t)\widetilde{N_K}^{-p}+t\widetilde{N_L}^{-p}}.
		\end{align*}
		It easy to check that 
		\begin{align*}
			(1-s)\cdot \widetilde{K}+_p s\cdot \widetilde{L}=\frac{1}{\left({(1-t)\widetilde{N_K}^{-p}+t\widetilde{N_L}^{-p}}\right)^{\frac{1}{p}}}\((1-t)\cdot K+_pt\cdot L\).
		\end{align*}
		From \eqref{equ:stab L2-BM-Wq=Wqp!0-q=0} for case $p\neq 0$, we have
			\begin{align}\label{equ:pf-1}
			&\frac{1}{p}\left(e^{\frac{p\left(\widetilde{E}((1-s)\cdot \widetilde{K} +_p s\cdot \widetilde{L})-\widetilde{E}(\widetilde{K})\right)}{|\mathbb{S}^n|}}-1\right)
			\nonumber
			\\
			\geq& s(1-s)|\mathbb{S}^n|^{-1}\left(\frac{2n+2+p}{2}-{\eta}\right)\delta_2^2( \widetilde{K},\widetilde{L}).
		\end{align}	
		Besides, we calculate that
		\begin{align}
			\label{equ:pf-2}
			&\widetilde{E}\((1-s)\cdot\tilde{ K} +_p s\cdot \tilde{L}\)-\widetilde{E}(\tilde{K})\nonumber
			\\
			=&\widetilde{E}\((1-t)\cdot K+_pt\cdot L\)-\frac{1}{p}|\mathbb{S}^n| \cdot 	\log\left((1-t)e^{p\frac{\widetilde{E}(K)}{|\mathbb{S}^n|}}+te^{p\frac{\widetilde{E}(L)}{|\mathbb{S}^n|}}\right).
		\end{align} 	
		Putting \eqref{equ:pf-2} into \eqref{equ:pf-1}, we obtain 
		\begin{align*}
			&\frac{1}{p}\left(
			e^{p\frac{\widetilde{E}\left((1-t)\cdot K+_pt\cdot L\right)}{|\mathbb{S}^n|}}
			-
			\left((1-t)e^{p\frac{\widetilde{E}(K)}{|\mathbb{S}^n|}}+te^{p\frac{\widetilde{E}(L)}{|\mathbb{S}^n|}}\right)
			\right)
			\\
			\geq& \frac{t(1-t)e^{\frac{p\widetilde{E}(K)}{|\mathbb{S}^n|}}e^{\frac{p\widetilde{E}(L)}{|\mathbb{S}^n|}} }{\left((1-t)e^{p\frac{\widetilde{E}(K)}{|\mathbb{S}^n|}}+te^{p\frac{\widetilde{E}(L)}{|\mathbb{S}^n|}}\right)\cdot|\mathbb{S}^n|}
			\left(\frac{2n+2+p}{2}-\eta\right)\delta_2^2(\tilde{K},\tilde{L}).
		\end{align*}
		That is \eqref{equ:stab L2-BM-q=0}.
		
		For $p=q=0$, we take $\tilde{K}=\widetilde{N_K}  K,\tilde{L}=\widetilde{N_L}  L$, where $\widetilde{N_K}=e^{-\frac{\widetilde{E}(K)}{|\mathbb{S}^n|}},  \widetilde{N_L}=e^{-\frac{\widetilde{E}(L)}{|\mathbb{S}^n|}}$. Then we know $\widetilde{E}(\tilde{K})=\widetilde{E}(\tilde{L})=\widetilde{E}(\mathbb{B}^{n+1})=0$ and $\tilde{K},\tilde{L}$ are also  $C^3$  convex bodies in $\mathcal{N}_{O(\varepsilon)}$ by Lemma \ref{lemma:C2-1-p-sum-C2-1}. 
		We can check that $\tilde{K}$ and $\tilde{L}$ satisfy $S_p(\tilde{K})=S_p(\tilde{L})={o}$. 
		Take
		\begin{align*}
			\gamma=\widetilde{N_K}^{1-t}\widetilde{N_L}^t,
		\end{align*} 
		then $(1-t)\cdot \tilde{K}+_0t\cdot \tilde{L}=\gamma \left((1-t)\cdot K+_0t\cdot L\right)$  by Lemma \ref{lemma:C2-1-p-sum-C2-1} (iii). 
		From \eqref{equ:stab L2-BM-Wq=Wq-q=0}, we have
		\begin{align}\label{equ:pf-3}
			&\widetilde{E}\((1-t)\cdot\tilde{ K} +_0 t\cdot \tilde{L}\)-\widetilde{E}(\tilde{K})\nonumber
			\\
			\geq&t(1-t)\left((n+1)-\eta\right)\delta_2^2(\tilde{K},\tilde{L}).
		\end{align}
		Since 
		\begin{align}\label{equ:pf-4}
			&\widetilde{E}\((1-t)\cdot\tilde{ K} +_0 t\cdot \tilde{L}\)=(\log\gamma)|\mathbb{S}^n|+\widetilde{E}\((1-t)\cdot K+_0 t\cdot L\)\nonumber
			\\
			=&\widetilde{E}\((1-t)\cdot K+_0 t\cdot L\)-\left((1-t)\widetilde{E}(K)+t\widetilde{E}(L)\right).
		\end{align}
		We substitute \eqref{equ:pf-4} into \eqref{equ:pf-3} and obtain \eqref{equ:stab L2-BM-p=q=0}.
	\end{proof}

	\section{A  uniqueness result of ($p,q$)-th dual curvature measure}\label{sec 6}
	The stability results for the inequalities imply the corresponding inequalities themselves. Using the dual $L^p$-Minkowski (or dual $L^p$-Brunn--Minkowski) inequalities for dual quermassintegrals and, when $q=0$, for the dual entropy, we then derive uniqueness result of ($p,q$)-th dual curvature measure.
	
	 Given that $\delta_2(K, NL)\geq 0$, moreover $\delta_2(K, NL) = 0$ if and only if $K = NL$. That is $K=\left(\frac{\widetilde{W}_q(K)}{\widetilde{W}_q(L)}\right)^{\frac{1}{q}}L$, which implies $K=cL$ for some constant $c>0$ (for $q\neq 0$). For $q=0$, we also have  $\delta_2(K, N_0L)\geq 0$, and  $\delta_2(K, N_0L) = 0$ if and only if $K = \exp\!\left((\widetilde E(K)-\widetilde E(L))
	 /|\mathbb S^n|\right)L$, which implies $K=cL$ for some constant $c>0$ (for $q=0$).
	Then Corollary \ref{thm:dual p-M-L2} yields the following  dual $L^p$-Minkowski inequality for $p > q-2n-2$.
	
	\begin{theorem}\label{thm:p-M-inequ}
		Let $p,q\in\mathbb{R}$ satisfy $2n+2+p-q>0$. There exists $\varepsilon>0$, that depends only on $n,p$ and $q$, such that if $C^3$ convex bodies $K,L\in\mathcal{N}_{\varepsilon}$ and $K,L$ satisfy  $S_p(K)=S_p(L)={o}$, then for  $q\neq 0, p\neq 0$, we have
		\begin{align}
			\label{equ:p-M-inequ-p!=0}
			\frac{1}{p}\left(\int_{\mathbb{S}^n}u_L^pu_K^{1-p}\rho_{K}^{q-n-1}\sigma_n(\lambda_K)~dx-\widetilde{W}_q(K)^{\frac{q-p}{q}}\widetilde{W}_q(L)^{\frac{p}{q}}\right)
			\geq
			0,
		\end{align}
		and for $q\neq 0, p=0$ we have
		\begin{align}
			\label{equ:p-M-inequ-p=0}
			\int_{\mathbb{S}^n}\left(\log(\frac{u_L}{u_K})\right)u_K\rho_{K}^{q-n-1}\sigma_n(\lambda_K) ~dx
			\geq
			\frac{1}{q}\widetilde{W}_q(K)\log\frac{\widetilde{W}_q(L)}{\widetilde{W}_q(K)}
			,
		\end{align}
		and for $p\neq 0, q=0$, 
		we have
		\begin{align}\label{equ:dual p-M-q=0}
			\frac{1}{p}\left(\int_{\mathbb{S}^n}u_L^pu_K^{1-p}\rho_K^{-n-1}\sigma_n(\lambda_K)~dx-|\mathbb{S}^n|e^{\frac{p\left(\widetilde{E}(L)-\widetilde{E}(K)\right)}{|\mathbb{S}^n|}}\right)
			\geq 0,
		\end{align}
		and for $p=0,q=0$, we have 
		\begin{align}
			\label{equ:dual p-M-p=0-q=0}
			\int_{\mathbb{S}^n}\log\left(\frac{u_L}{u_K}\right)u_K\rho_K^{-n-1}\sigma_n(\lambda_K)~dx+\widetilde{E}(K)-\widetilde{E}(L)
			\geq 0.
		\end{align}
		Moreover, equality holds in \eqref{equ:p-M-inequ-p!=0}$\thicksim$\eqref{equ:dual p-M-p=0-q=0} if and only if $K=cL$ for some $c>0$. 
	\end{theorem}
	\begin{proof}[\textbf{Proof of Corollary \ref{thm:unique}}]
		(1).  When $p\neq 0, q\neq 0$, by \eqref{equ:Christolell-equ}, we can derive
		\begin{align}\label{equ:pf-christ-0}
			\int_{\mathbb{S}^n}u_L^pu_K^{1-p}\rho_{K}^{q-n-1}\sigma_n(\lambda_K)~dx=\int_{\mathbb{S}^n}u_L^pu_L^{1-p}\rho_{L}^{q-n-1}\sigma_n(\lambda_L)~dx=\widetilde{W}_q(L).
		\end{align}
		Substituting \eqref{equ:pf-christ-0} into \eqref{equ:p-M-inequ-p!=0}, we obtain an equality when $p= q$, and the equality condition gives $K=cL$. 
		
		If $p\neq q$, putting \eqref{equ:pf-christ-0} into \eqref{equ:p-M-inequ-p!=0}, we have
		\begin{align}\label{equ:pf-christ-1}
			\frac{1}{p}\widetilde{W}_q(L)\left(1-\left(\frac{\widetilde{W}_q(K)}{\widetilde{W}_q(L)}\right)^{\frac{q-p}{q}}\right)\geq 0.
		\end{align}
		The same conclusion holds upon interchanging $K$ and $L$  in Theorem \ref{thm:p-M-inequ}, analogously, we have
		\begin{align}
			\label{equ:pf-christ-2}		\frac{1}{p}\widetilde{W}_q(K)\left(1-\left(\frac{\widetilde{W}_q(L)}{\widetilde{W}_q(K)}\right)^{\frac{q-p}{q}}\right)\geq 0.
		\end{align}
		Combining \eqref{equ:pf-christ-1} and \eqref{equ:pf-christ-2}, we obtain
		\begin{align*}
			\widetilde{W}_q(K)=\widetilde{W}_q(L).
		\end{align*}
		Then 
		\begin{align}\label{equ:pf-christ-3}
			\widetilde{W}_q(K)^{\frac{q-p}{q}}\widetilde{W}_q(L)^{\frac{p}{q}}	=\widetilde{W}_q(L).
		\end{align}
		Observing \eqref{equ:pf-christ-0} and \eqref{equ:pf-christ-3} shows that equality holds in \eqref{equ:p-M-inequ-p!=0}, whence Theorem \ref{thm:p-M-inequ} implies  $K=cL$,  putting it into \eqref{equ:Christolell-equ} yields  $c=1$ for $p\neq q$.  
		\\
		(2). When $p= 0, q\neq 0$,  from Theorem \ref{thm:p-M-inequ}, 
		we obtain
		\begin{align}
			\int_{\mathbb{S}^n} \Bigl(\log \frac{u_L}{u_K}\Bigr) u_K \rho_{K}^{q-n-1}\sigma_n(\lambda_K) ~dx
			- \frac{1}{q}\widetilde{W}_q(K) \log\frac{\widetilde{W}_q(L)}{\widetilde{W}_q(K)} \geq 0. \label{equ:pf-Chri-2}
		\end{align}
		Interchanging $K$ and $L$ in Theorem \ref{thm:p-M-inequ} gives
		\begin{align*}
			\int_{\mathbb{S}^n} \Bigl(\log \frac{u_K}{u_L}\Bigr) u_L \rho_{L}^{q-n-1}\sigma_n(\lambda_L) ~ dx
			- \frac{1}{q}\widetilde{W}_q(L) \log\frac{\widetilde{W}_q(K)}{\widetilde{W}_q(L)} \geq  0. 
		\end{align*}
		Using the condition \eqref{equ:Christolell-equ} for $p=0$, we derive that
		\begin{align}
			-\int_{\mathbb{S}^n} \Bigl(\log \frac{u_L}{u_K}\Bigr) u_K \rho_{K}^{q-n-1}\sigma_n(\lambda_K) ~dx
			- \frac{1}{q}\widetilde{W}_q(L) \log\frac{\widetilde{W}_q(K)}{\widetilde{W}_q(L)} \geq  0. \label{equ:pf-Chri-4}
		\end{align}
		and
		\begin{align}
			\widetilde{W}_q(K)=\int_{\mathbb{S}^n}u_K\rho_{K}^{q-n-1}\sigma_n(\lambda_K)~dx\overset{\eqref{equ:Christolell-equ}}{=}\int_{\mathbb{S}^n}u_L\rho_{L}^{q-n-1}\sigma_n(\lambda_L)~dx= \widetilde{W}_q(L). \label{equ:pf-Chri-7}
		\end{align}
		Substituting \eqref{equ:pf-Chri-7} into \eqref{equ:pf-Chri-2} gives
		\begin{align}
			\int_{\mathbb{S}^n} \Bigl(\log \frac{u_L}{u_K}\Bigr) u_K \rho_{K}^{q-n-1}\sigma_n(\lambda_K) ~dx \geq  0, \label{equ:pf-Chri-8}
		\end{align}
		while substituting it into \eqref{equ:pf-Chri-4} gives
		\begin{align*}
			-\int_{\mathbb{S}^n} \Bigl(\log \frac{u_L}{u_K}\Bigr) u_K \rho_{K}^{q-n-1}\sigma_n(\lambda_K) ~dx \geq  0. 
		\end{align*}
		Hence the integral must be zero, and equality holds in \eqref{equ:pf-Chri-2}.
		
		By Theorem \ref{thm:p-M-inequ}, equality in \eqref{equ:pf-Chri-2} occurs if and only if $K = cL$, putting it into \eqref{equ:Christolell-equ} yields have $c=1$.  
		\\
		(3).  When $p \neq 0$ and $q = 0$, substituting the condition \eqref{equ:Christolell-equ} (for $q = 0$) into \eqref{equ:dual p-M-q=0} and using $\int_{\mathbb{S}^n} u_L \rho_L^{-n-1} \sigma_n(\lambda_L) ~ dx = |\mathbb{S}^n|$ yields
		\begin{align*}
			\widetilde{E}(K)-\widetilde{E}(L)
			\geq 0.
		\end{align*}
		The same conclusion holds if we interchange $K$ and $L$ in Theorem \ref{thm:p-M-inequ}, analogously, we obtain
		\begin{align*}
			\widetilde{E}(L)-\widetilde{E}(K)
			\geq 0.
		\end{align*}
		Then $\widetilde{E}(L) = \widetilde{E}(K)$, which further implies that equality holds in \eqref{equ:dual p-M-q=0}, and consequently $K = cL$, putting it into \eqref{equ:Christolell-equ} yields  $c=1$ for $p\neq0= q$.  
		\\
		(4). When $p=q=0$, putting the condition \eqref{equ:Christolell-equ} (for $p=q = 0$) into \eqref{equ:dual p-M-p=0-q=0}  yields
		\begin{align*}
			\widetilde{E}(L)-\widetilde{E}(K)	\overset{\eqref{equ:dual p-M-p=0-q=0}}{\leq}& \int_{\mathbb{S}^n}\log\left(\frac{u_L}{u_K}\right)u_K\rho_K^{-n-1}\sigma_n(\lambda_K)~dx
			\\
			\overset{\eqref{equ:Christolell-equ}}{=}&-\int_{\mathbb{S}^n}\log\left(\frac{u_K}{u_L}\right)u_L\rho_L^{-n-1}\sigma_n(\lambda_L)~dx
			\\
			\overset{\eqref{equ:dual p-M-p=0-q=0}}{\leq} &- \left(\widetilde{E}(K)-\widetilde{E}(L)\right)=\widetilde{E}(L)-\widetilde{E}(K).
		\end{align*}
		Thus equality holds in \eqref{equ:dual p-M-p=0-q=0},
		so $K=cL$ for some $c>0$. 
	\end{proof}

		\vskip 2mm  {\bf Acknowledgment.} 
	The research of  JG was supported by NSF of Zhejiang province No. LQN26A010008.
	The research of  GL was supported by NSFC No. 12261105.
	
	\vskip 2mm
	\textbf{Data Availability.} Data sharing is not applicable to this article as no datasets were generated or analysed during the current study.
	
	\vskip 2mm
	\textbf{Conflict of interest.} On behalf of all authors, the corresponding author states that there is no conflict of interest.


	\section*{Reference}
	

	\bigbreak
	
	\noindent Shanwei Ding,  {\tt dingsw@whu.edu.cn} \\
	School of Mathematics and Statistics, Huazhong University of Science and
	Technology, Wuhan, Hubei 430074, China\\

	\noindent  Jinyu Gao, {\tt gaojy@zstu.edu.cn} \\
	Department of Mathematics,
	Zhejiang Sci-Tech University,
	Hangzhou, Zhejiang 310018, China\\

	\noindent  Guanghan Li, {\tt ghli@whu.edu.cn} \\
	School of Mathematics and Statistics, Wuhan University, 
	Wuhan, Hubei 430072, China

\end{document}